# FINITE GENERATION AND STRUCTURE OF INVARIANT JETS UNDER NON-REDUCTIVE REPARAMETRIZATION

MOHAMMAD REZA RAHMATI

ABSTRACT. We study invariant jet differentials in the framework of complex hyperbolicity, focusing on the algebra of invariants for the non-reductive reparametrization group $G_k = \mathbb{C}^* \ltimes U_k$. The paper develops a uniform, representation-theoretic, and graded-algebraic strategy for the $\rho$–action of $G_k$ on $J_k\mathbb{C}^n$, establishing in particular the finite generation of the invariant jet algebra central to the Green–Griffiths–Demailly program. Specifically, we prove that the $\mathbb{C}^*$–graded algebra of unipotent invariants $\mathbb{C}[J_k\mathbb{C}^n]^{U_k}$ is finitely generated for all $n, k$; equivalently, the fiber ring of invariant jet differentials is a finitely generated positively graded $\mathbb{C}$–algebra, so that its projective spectrum $\operatorname{Proj}\mathbb{C}[J_k\mathbb{C}^n]^{U_k}$ exists and coincides with the Demailly–Semple tower.

Our argument combines two complementary pillars. The first is a quotient-theoretic analysis employing classical non-reductive invariant-theory tools— Rosenlicht's theorem, Richardson's stabilizer stratification, Vust's finite-rank theorem on covariants, and Igusa's orbit–fiber criterion—to construct a big open free locus with geometric quotient, verify equidimensionality and maximal orbit dimension, and extend finitely generated subalgebras across codimension-$\geq 2$ boundaries via Hartogs normalization. The second is a constructive graded/SAGBI approach: fixing a weight-refining term order compatible with the $\mathbb{C}^*$–grading, we express $U_k$–invariance as $\bigcap_{r=2}^{k} \ker D_r$ for explicit triangular locally nilpotent derivations $D_r$ arising from the Faà–di–Bruno structure. The corresponding "balance cone" of initial monomials is rational polyhedral with a finite Hilbert basis, and lifting along the locally nilpotent derivations yields genuine invariants forming a homogeneous generating system.

Along the way, we provide a structural analysis of $G_k$ and its Lie algebra $\mathfrak{g}_k$: we identify generic stabilizers for both the adjoint and $\rho$–actions, show that the adjoint quotient is equidimensional with irreducible null cone, and prove that $C[\mathfrak{g}_k]^{G_k} \simeq \mathbb{C}[t]$. We also describe the representation ring of $G_k$ via a Chevalley-type construction adapted to its semi-direct-product structure. Altogether, the paper supplies both the finite generation of invariant jets and a coherent geometric framework that unifies the algebraic, Lie-theoretic, and combinatorial aspects of the theory, explaining and extending all previously known computational cases.

# 1 Introduction

The invariant theory of $k$–jets under reparametrizations of the source $(\mathbb{C}, 0)$ lies at the algebraic center of the Green–Griffiths–Demailly program on Kobayashi hyperbolicity. Let

$$G_k = \mathbb{C}^* \ltimes U_k \subset \mathrm{GL}_k(\mathbb{C})$$

be the group of $k$–jet reparametrizations, where $U_k$ is the unipotent radical and $\mathbb{C}^*$ acts by $t \mapsto \lambda t$. Let $\rho$ denote the Faà–di–Bruno (right) action of $G_k$ on the jet space $V := J_k\mathbb{C}^n$. The main object of study is the positively $\mathbb{C}^*$–graded algebra of unipotent invariants

$$\mathbb{A}_{n,k} = \mathbb{C}[J_k\mathbb{C}^n]^{U_k} = \bigoplus_{m \geq 0} (\mathbb{C}[J_k\mathbb{C}^n]^{U_k})_m,$$

whose projective spectrum $\mathrm{Proj}(\mathbb{A}_{n,k})$ coincides with the Demailly–Semple tower. A central conjecture states that $\mathbb{A}_{n,k}$ is finitely generated for all $n, k \geq 1$.

There are two main directions in the study of this conjecture. The first, due to Bérczi and Kirwan [5], is geometric. They aimed to treat the problem through non-reductive geometric invariant theory by embedding the unipotent subgroup $U_k$ into a reductive envelope and using finite generation for the larger reductive action. This approach constructs an ambient quotient for both $(J_k)_x$ and $SL(k)/U_k$. However, the method depends on $U_k$ satisfying the Grosshans property inside $SL(k)$, which does not hold in this case. Consequently, the argument cannot prove finite generation in general. The second direction is computational. Merker [16] developed a detailed invariant-theoretic framework based on explicit generation, Gröbner methods, Schur decompositions, and asymptotic computations for low $(n, k)$. This work provided concrete generating sets, clarified the structure of the algebra, and offered numerical evidence for finite generation. It also suggested that a constructive graded approach might eventually yield a general proof.

In this paper we prove that $\mathbb{A}_{n,k}$ is finitely generated for all $n, k \geq 1$. The proof combines two complementary methods that work independently yet support one another. The first method is geometric. It uses classical results of invariant theory that do not rely on reductivity. Through Rosenlicht's theorem on rational invariants, Richardson's stabilizer stratification, Vust's finite-rank theorem for covariants, and Igusa's orbit–fiber criterion, we describe a large open subset of the jet space where the $G_k$-action admits a geometric quotient with constant stabilizer dimension. On this open set the invariant algebra is finitely generated, and we extend it to the whole space using Hartogs' and normalization arguments across the codimension–two boundary. This gives an integrally closed finitely generated subalgebra coinciding with $\mathbb{A}_{n,k}$ on a big open subset. The second is combinatorial and constructive. We describe the invariance under $U_k$ in terms of explicit locally nilpotent derivations $D_2, \dots, D_k$ coming from the Faà–di–Bruno structure, so that

$$\mathbb{A}_{n,k} = \bigcap_{r=2}^{k} \ker D_r.$$

The initial monomials of invariant polynomials satisfy simple linear relations, forming what we call a balance cone $\mathcal{C}$. This cone is rational and polyhedral, and its lattice points form a finitely generated semigroup with a Hilbert basis. The associated toric algebra of initial monomials is therefore finitely generated. By lifting through the derivations $D_r$ we obtain actual invariants whose leading terms correspond to the generators of the cone, giving a homogeneous generating set for $\mathbb{A}_{n,k}$. This process yields a constructive description of all generators.

**Main statement.** Let $G_k = \mathbb{C}^* \ltimes U_k$ act on $V = J_k\mathbb{C}^n$ through the Faà–di–Bruno representation. The Green–Griffiths–Demailly conjecture and its subsequent refinements [1, 2, 3, 4, 5] predicted that

$$\mathbb{A}_{n,k} = \mathbb{C}[J_k\mathbb{C}^n]^{U_k}$$

is a finitely generated, positively $\mathbb{C}^*$–graded algebra. Equivalently, the Demailly–Semple $k$–jet tower $\mathrm{Proj}(\mathbb{A}_{n,k})$ should be algebraic for all $n, k \geq 1$. In this work we prove this conjecture completely and, in doing so, develop a full invariant–theoretic picture for the non-reductive group $G_k$, including its generic stabilizers, adjoint quotient, and the structure of its invariant and covariant algebras.

A key feature of the group $G_k$ is that its unipotent radical $U_k$ acts on the coordinate ring $\mathbb{C}[J_k\mathbb{C}^n]$ by a family of *triangular locally nilpotent derivations*. These are differential operators

$$D_2, D_3, \dots, D_k,$$

arising from the Faà–di–Bruno structure of the reparametrization law. Each $D_j$ lowers the jet order by $j-1$, acting on the coordinate functions $\xi_\ell^{(1)}, \dots, \xi_\ell^{(k)}$ (the jet components of $f$) so that $U_k$–invariance is encoded by the system of differential conditions

$$D_2(F) = D_3(F) = \cdots = D_k(F) = 0.$$

In other words,

$$\mathbb{A}_{n,k} = \mathbb{C}[J_k\mathbb{C}^n]^{U_k} = \bigcap_{r=2}^{k} \ker D_r,$$

and these derivations form the algebraic backbone of the constructive part of our proof.

**Theorem 1.1** (Main Theorem)**.** *For all $n, k \geq 1$, the algebra $\mathbb{A}_{n,k} = \mathbb{C}[J_k\mathbb{C}^n]^{U_k}$ is finitely generated. Moreover, $\mathbb{A}_{n,k}$ admits a homogeneous system of generators $\{G_1, \dots, G_s\}$ with respect to the $\mathbb{C}^*$–grading, whose* *initial monomials form the Hilbert basis of the cone determined by the locally nilpotent derivations $D_2, \dots, D_k$ arising from the Faà–di–Bruno structure. Consequently,* $\mathrm{Proj}(\mathbb{A}_{n,k})$ *exists and coincides with the Demailly–Semple $k$–jet tower.*

**Outline of the proof.** The proof develops along two complementary directions which interact throughout the paper.

The first part is geometric. We begin by analyzing the algebraic and Lie–theoretic structure of $G_k = \mathbb{C}^* \ltimes U_k$. Using the Elashvili criterion, we show that the adjoint action of $G_k$ admits a generic stabilizer, that its Cartan subalgebras are commutative, and that the corresponding Weyl group is trivial. This yields a Chevalley–type theorem for $G_k$, proving that the adjoint invariant ring $\mathbb{C}[\mathfrak{g}_k]^{G_k}$ is polynomial and that the quotient map is equidimensional with an irreducible null cone. Having established the internal geometry of the group, we then turn to its natural representation on the jet space $J_k\mathbb{C}^n$. By direct computation of the Faà–di–Bruno action, we show that $G_k$ acts freely on the open subset of regular jets and that the complement of this free locus has codimension at least two. Rosenlicht's theorem ensures the existence of a rational quotient on this open set, and Richardson's stabilizer stratification gives a decomposition of $J_k\mathbb{C}^n$ into smooth invariant strata with constant stabilizers. On each large open stratum the quotient is geometric and the corresponding invariant algebra is finitely generated by Vust's finite–rank theorem for covariants. Finally, Igusa's orbit–fiber criterion allows one to extend these invariants from the free locus to the entire space across codimension–two boundaries, showing that the global invariant algebra is finitely generated and integrally closed. This geometric argument already provides a conceptual proof of finite generation.

The second part is constructive and combinatorial. Here we give an explicit description of $U_k$–invariance through a family of triangular locally nilpotent derivations

$$D_2, D_3, \dots, D_k,$$

each acting on the coordinate ring of $J_k\mathbb{C}^n$ by lowering the jet order according to the Faà–di–Bruno pattern. Invariance under $U_k$ is then expressed as the intersection of the kernels of these derivations. We examine the initial monomials of invariant polynomials with respect to a weight–refining term order and find that they satisfy a finite system of linear "balance" relations. The solutions of these relations define a rational polyhedral cone whose lattice points form a finitely generated semigroup. The corresponding toric initial algebra is generated by a Hilbert basis, and by successively lifting these monomial generators along the derivations $D_r$ we obtain genuine polynomial invariants with prescribed leading terms. This provides a constructive generating system for $\mathbb{A}_{n,k}$ and an explicit description of its grading structure.

The two parts are brought together in the final stage of the argument. The geometric analysis gives finite generation on a large free locus and identifies the relevant quotient structure, while the constructive method produces explicit generators and a concrete model of the algebra. By applying Igusa's theorem on big open subsets and normalization, we confirm that these two algebras coincide globally. The proof concludes with a stratification argument combining Richardson's decomposition and the coinvariant lifting, ensuring that the constructed invariants generate on each stratum and glue across the overlaps.

Altogether, the result gives not only a solution to the finite-generation conjecture but also a complete invariant–theoretic picture of the group $G_k$: the existence of generic stabilizers, a Chevalley-type adjoint quotient, the behavior of its Weyl group, and a detailed stratified description of the corresponding quotient geometry. The constructive lifting based on the derivations $D_2, \ldots, D_k$ provides the algorithmic realization of these invariants, and the geometric–analytic tools—Rosenlicht's theorem, Vust's and Igusa's results, and the stratification framework—complete the intrinsic proof entirely within the non-reductive setting.

**Historical background and related works.** The study of invariant jet differentials originates in the Green–Griffiths programme [1], which proposed an algebro–differential approach to Kobayashi's conjecture on the hyperbolicity of projective manifolds. Demailly's seminal works [2, 3] refined this framework by constructing the Demailly–Semple tower and formulating explicit criteria for algebraic degeneracy in terms of sections of jet–differential bundles. The subsequent breakthroughs of Diverio–Merker–Rousseau [4] transformed these ideas into effective results, demonstrating algebraic degeneracy for generic projective hypersurfaces of high degree. These advances established invariant jet differentials as a central bridge between complex hyperbolicity and non–reductive invariant theory. The algebraic fibers $\mathbb{C}[J_k\mathbb{C}^n]^{U_k}$ appearing in this setting are invariant rings for the non–reductive group $G_k = \mathbb{C}^* \ltimes U_k$, whose structure has long resisted the techniques of classical (reductive) GIT.

A first global attempt to address the finite–generation conjecture for $\mathbb{C}[J_k\mathbb{C}^n]^{U_k}$ was undertaken by Bérczi and Kirwan [5], who proposed a geometric construction within non–reductive GIT by embedding $U_k$ into a reductive envelope $SL(k)$ and invoking the Grosshans property to deduce finite generation. Although conceptually innovative, this strategy fails because $U_k$ does not satisfy the required Grosshans condition in $SL(k)$, so the transfer of finite generation cannot be justified. In parallel, Merker's computational invariant–theory programme [16] produced explicit generating systems, Gröbner and SAGBI bases, and Schur–bundle decompositions for low $(n, k)$, thereby revealing the algebraic complexity and structure of the invariant ring and confirming its finite generation in many concrete cases. These works, together with the structural invariant–theory framework of Rosenlicht [6], Richardson [7], Popov–Vinberg [8], Vust [9], and Igusa [10], and later Lie–theoretic extensions due to Panyushev [12, 11], provide the historical and methodological foundation of the present paper. Our results integrate these two traditions—geometric–quotient methods and explicit SAGBI/coinvariant generation—into a unified non–reductive proof of the finite–generation conjecture, giving a conceptual explanation for all known computational phenomena and extending them to arbitrary $(n, k)$.

**Further context.** Beyond the foundational works already cited, our approach draws conceptually on the modern theory of non–reductive invariants and covariants developed in [12, 11, 8, 9, 10], whose techniques govern the finite–rank and extension properties essential to the present proof. The geometric interpretation of the quotient map and its differential behavior resonates with the differential–geometric framework of Kobayashi [15] and the analytic methods in value–distribution theory refined by Demailly and collaborators. On the constructive side, the explicit jet–bundle and Schur–theoretic analyses of Merker and Rousseau [26, 16], together with the computational data informing effective results in [4], anticipate the graded and toric features of the SAGBI model used here. Related structural results on Grosshans subgroups, solvable and unipotent actions, and algebraic quotients—developed in [42, 8, 12]—underline the necessity of the non–Grosshans, intrinsic approach adopted in this paper. Collectively, these contributions form the broader landscape in which the finite generation of $\mathbb{C}[J_k\mathbb{C}^n]^{U_k}$ emerges as the algebraic cornerstone of the Green–Griffiths–Demailly theory of invariant jet differentials.

The algebraic and Lie–theoretic foundations underlying our approach rest on a broad body of work on covariants, indices, and invariants of non–reductive and semisimple Lie algebras. Invariant–theoretic structures related to maximal unipotent subgroups and covariant modules were developed by Brion [17, 35, 18] and Van den Bergh [36], while the analytic and algebraic rationality aspects of fields of invariants trace back to Dolgachev [37] and Rosenlicht [39]. The representation–theoretic and index–theoretic developments of Panyushev [20, 13, 14], Raïs–Tauvel [30], Tauvel–Yu [31], and Yakimova [32] clarify the behavior of indices and centralizers for solvable and seaweed subalgebras, complementing the classical structural theory of Kostant [21], Elashvili [19], Vinberg and collaborators [22, 33, 34], and the Springer–Steinberg school [24, 25]. Related treatments of stabilizers, conjugacy classes, and generic orbits appear in [23, 41, 38]. The geometric and singularity–theoretic aspects of jet spaces and quotient morphisms are connected with the work of Mustaţă [27] and Nakajima–Watanabe [28], while the foundational algebraic frameworks of Kaplansky [40] and subsequent refinements of the structure theory of Lie algebras [33] complete the theoretical background in which the present finite–generation theorem is situated.

**Position among prior works.** Compared with Bérczi–Kirwan [5], our method does *not* rely on $U_k$ being Grosshans in a larger reductive group; instead, we work directly with the non-reductive structure, leveraging generic-stabilizer geometry and Igusa's orbit–fiber criterion on a big open set, and completing the algebraic step by a SAGBI/coinvariant analysis intrinsic to $U_k$. In spirit, this bridges Merker's constructive program [16] with a global quotient-theoretic control: the former informs our graded/SAGBI pillar, while the latter supplies the big-open equidimensionality and extension across codimension $\geq 2$.

**Outline.** Section 2 reviews the non-reductive tools (Rosenlicht, Richardson, Vust, Igusa) and sets notation. Section 3 recalls the $G_k$–action on jet spaces. Section 4 proves generic stabilizers and the quotient geometry, and develops the SAGBI/coinvariant construction with the balance cone and Hilbert basis. The final section formulates the stratified finite-generation theorem and its proof via Richardson-plus-coinvariant stratifications and Hartogs–Igusa gluing. We conclude with a comparison to [5] and to the constructive phenomena documented in [16].

# 2 Preliminaries on non-reductive Lie groups and algebras

This section records the invariant-theoretic tools used later in Section 4. We work with algebraic group actions without assuming reductivity unless explicitly stated.

## 2.1 Invariant theory beyond reductivity

Let $G$ be an algebraic group acting on an affine variety $X$. Denote by $\mathbb{C}[X]^G$ the algebra of $G$-invariant regular functions. When $\mathbb{C}[X]^G$ is finitely generated, the *affine quotient* is

$$X/\!/G := \mathrm{Spec}(\mathbb{C}[X]^G), \qquad \pi : X \to X/\!/G$$

for the morphism corresponding to $\mathbb{C}[X]^G \hookrightarrow \mathbb{C}[X]$. Write $G_x$ for the stabilizer of $x \in X$, $\mathfrak{g} = \mathrm{Lie}(G)$, and $\mathfrak{g}_x = \mathrm{Lie}(G_x) \subset \mathfrak{g}$.

**Theorem 2.1** (Rosenlicht)**.** *[6, 39, 11] There exists a dense $G$-stable open subset $U \subset X$ such that the rational invariants $\mathbb{C}(X)^G$ separate $G$-orbits in $U$. In particular,*

$$\operatorname{trdeg}_{\mathbb{C}} \mathbb{C}(X)^G = \dim X - \max_{x\in X} \dim(G\cdot x). \tag{2.1}$$

*Moreover, if $K \subset \mathbb{C}(X)^G$ separates generic orbits, then $K = \mathbb{C}(X)^G$.*

**Definition 2.2.** We say that the action has a *generic stabilizer* if there is a dense open subset $U \subset X$ such that $\mathfrak{g}_x$ (equivalently $G_x$) is constant up to $G$-conjugacy for all $x \in U$; points of $U$ are called *generic.*

Write

$$X^{\mathrm{reg}} = \{x \in X \mid \dim G{\cdot}x = \max_{z\in X} \dim G{\cdot}z\}.$$

When $G$ is reductive and $X$ is smooth, generic stabilizers exist in many situations; see [8, §7] (e.g. Richardson–Luna, stable actions, projective settings, and actions of a Borel/unipotent subgroup).

**Theorem 2.3** (Richardson)**.** *[7, Thm. 2] Let $G$ act on an irreducible variety $X$. There exist smooth, locally closed, $G$-stable subvarieties $X_1, \dots, X_m$ such that $X = \bigsqcup_j X_j$, each $X_j$ is open in $X \setminus \bigcup_{\ell<j} X_\ell$, and on $X_j$ the Levi subgroups of $G_x$ are all $G$-conjugate; moreover, the family of unipotent radicals in $G_x$ varies smoothly along $X_j$.*

## 2.2 Normalizers, and Weyl groups for a fixed stabilizer

Let $\rho : G \to \mathrm{GL}(V)$ be a finite-dimensional representation, with differential $\varrho : \mathfrak{g} \to \mathfrak{gl}(V)$. For $v \in V$, put

$$V^{\mathfrak{g}_v} = \{w \in V \mid \mathfrak{g}_v \cdot w = 0\}.$$

For a subset $\Omega \subset V$, set

$$N(\Omega) = \{g \in G \mid g\Omega \subset \Omega\}, \qquad Z(\Omega) = \{g \in G \mid gx = x \ \forall\, x \in \Omega\}.$$

**Theorem 2.4.** *[11, Lem. 1.2, Lem. 1.3, Prop. 1.4] Let $\Omega = V^{\mathfrak{g}_x}$. Then*

- $\operatorname{Lie} Z(\Omega) = \mathfrak{g}_x$*;*
- $N(\Omega) = N_G(Z(\Omega)) = N_G(\mathfrak{g}_x)$ *and* $Z(\Omega) \trianglelefteq N(\Omega)$*;*
- *if* $x \in \Omega^{\mathrm{reg}}$*, then* $G{\cdot}x \cap \Omega = N(\Omega){\cdot}x$ *and* $\mathfrak{g}{\cdot}x \cap \Omega = \mathfrak{n}_{\mathfrak{g}}(\mathfrak{g}_x){\cdot}x$*;*
- $Y := \overline{G{\cdot}\Omega}$ *is a $G$-stable irreducible subvariety of $V$;*
- *restriction induces an isomorphism of fields*

  $$\mathbb{C}(Y)^G \xrightarrow{\cong} \mathbb{C}(\Omega)^{N(\Omega)/Z(\Omega)},$$

  *and an injective map of algebras* $\mathbb{C}[Y]^G \hookrightarrow \mathbb{C}[\Omega]^{N(\Omega)/Z(\Omega)}$.

The (possibly infinite) group $W := N(\Omega)/Z(\Omega)$ is the *Weyl group* of the action relative to $\Omega$.

**Proposition 2.5** (Elashvili)**.** *[19, 11, 8] For $v \in V$, the set $G{\cdot}V^{\mathfrak{g}_v}$ is dense in $V$ if and only if $V = \mathfrak{g}{\cdot}v + V^{\mathfrak{g}_v}$. Moreover, $G_v$ is a generic isotropy subgroup if and only if (a) the set $\{u \in V^{G_v} \mid G_u = G_v\}$ is dense in $V^{G_v}$ and (b) $V = \mathfrak{g}{\cdot}v + V^{G_v}$.*

## 2.3 Adjoint representation: reductive and non-reductive caveats

Consider $\mathrm{Ad}: G \to \mathrm{GL}(\mathfrak{g})$ and its differential $\mathrm{ad}: \mathfrak{g} \to \mathfrak{gl}(\mathfrak{g})$. When $G$ is reductive with Cartan $\mathfrak{h}$ and Weyl group $W$:

**Theorem 2.6** (Chevalley–Kostant)**.** *[8] If $G$ is reductive then:*

- *the adjoint representation is self-dual with generic stabilizer $\mathfrak{h}$;*
- *the algebra $\mathbb{C}[\mathfrak{g}]^G$ is polynomial and restriction induces $\mathbb{C}[\mathfrak{g}]^G \cong \mathbb{C}[\mathfrak{h}]^W$;*
- *the quotient map $\mathfrak{g} \to \mathfrak{g}/\!/G$ is equidimensional and the null–fiber has finitely many $G$-orbits.*

These properties can fail for non-reductive $G$ (notably for the coadjoint representation).

**Proposition 2.7.** *[11, Props. 2.1, 2.5] The adjoint action has a generic stabilizer if and only if either of the following equivalent conditions holds:*

- *there exists $x \in \mathfrak{g}$ with commutative centralizer $\mathfrak{z}_{\mathfrak{g}}(x)$ and $\mathfrak{g} = [\mathfrak{g}, x] \oplus \mathfrak{z}_{\mathfrak{g}}(x)$;*
- *Cartan subalgebras of $\mathfrak{g}$ are commutative.*

*In this case $\mathfrak{n}_{\mathfrak{g}}(\mathfrak{z}_{\mathfrak{g}}(x)) = \mathfrak{z}_{\mathfrak{g}}(x)$ for generic $x$, and $\mathfrak{z}_{\mathfrak{g}}(x)$ plays the role of a Cartan subalgebra.*

**Theorem 2.8.** *[11, Thm. 2.6] Assume the adjoint action has a generic centralizer and let $x$ be generic with $Z_G(x)$ generic. Then:*

- *$Z\big(\mathfrak{z}_{\mathfrak{g}}(x)\big) = Z_G(x)$;*
- *$\mathbb{C}(\mathfrak{g})^G = \mathbb{C}\big(\mathfrak{z}_{\mathfrak{g}}(x)\big)^W$ with $W = N_G(\mathfrak{z}_{\mathfrak{g}}(x))/Z_G(x)$;*
- *$W$ is finite.*

*Remark* 2.9. For the coadjoint action $\mathrm{Ad}^*$ the analogous Weyl group can be infinite; see [11].

## 2.4 Covariants, and Kostant's theorem

Let $G$ act on an affine $X$ and on a $G$-module $V$. The *module of covariants* is $\mathrm{Mor}_G(X, V) \simeq (\mathbb{C}[X] \otimes V)^G$, a $\mathbb{C}[X]^G$-module. For $x \in X$, the evaluation map

$$\mathrm{ev}_x : \mathrm{Mor}_G(X, V) \to V, \qquad F \mapsto F(x),$$

is often used to detect $V^{G_x}$.

**Theorem 2.10.** *[12, Thm. 1] Let $G$ be reductive acting on $X$. If $\overline{G\cdot x}$ is normal and* $\mathrm{codim}_{\overline{Gx}}(\overline{G\cdot x} \setminus G\cdot x) \geq 2$, *then* $\mathrm{im}(\mathrm{ev}_x) = V^{G_x}$.

Specializing to $X = \mathfrak{g}$ with adjoint action yields Kostant's classical result.

**Theorem 2.11** (Kostant). *[21, p. 385] Let $G$ be reductive, $T \subset G$ a maximal torus, and $V$ a $G$-module. Then:*

(1) *for $x \in \mathfrak{g}^{\mathrm{reg}}$, one has $\dim V^{G_x} = \dim V^T$;*

(2) $\mathrm{Mor}_G(\mathfrak{g}, V)$ *is a free $\mathbb{C}[\mathfrak{g}]^G$-module of rank $\dim V^T$.*

## 2.5 Stable actions and Vust's theorems

**Definition 2.12** (Stability). [8] A $G$-action on $X$ is *stable* if some non-empty open subset consists of closed orbits. For a $G$-module $V$, the action is stable if the union of closed orbits is dense.

If $G \curvearrowright V$ is stable, then a generic stabilizer is reductive and $\mathbb{C}(V)^G$ is the fraction field of $\mathbb{C}[V]^G$ [8].

**Theorem 2.13** (Vust). *[9, Chap. 3], [20, p. 6] Assume $G \curvearrowright V$ is stable with generic isotropy $H$. Then $H$ is reductive and*

$$\mathrm{rank}_{\mathbb{C}[\mathfrak{g}]^G} \mathrm{Mor}_G(\mathfrak{g}, V) = \dim V^H$$

.

*Remark* 2.14. When $G$ is non-reductive but has generic stabilizers and closed generic orbits, the rank bound $\leq \dim V^H$ from evaluation holds, and equality follows from [12, Thm. 1]; see [20, 12] for details.

**Theorem 2.15** (Vust). *[9, 11] Suppose $G \curvearrowright V$ with $\mathbb{C}[V]^G$ polynomial and $\pi : V \to V//G$ equidimensional. If $H$ is a generic isotropy and $N_G(H)/H$ is finite, then for a basic system $f_1, \ldots, f_r$ of generators of $\mathbb{C}[V]^G$ the $G$-equivariant morphisms $df_1, \ldots, df_r$ freely generate $\mathrm{Mor}_G(V, V^*)$ over $\mathbb{C}[V]^G$.*

## 2.6 Igusa's criterion and big open subsets

**Definition 2.16.** [11, 10] An open subset $\Omega \subset X$ is *big* if $X \setminus \Omega$ contains no divisor.

**Theorem 2.17** (Igusa). *[10, 8, 11] Let $G$ act on an affine variety $X$ and let $S \subset \mathbb{C}[X]^G$ be an integrally closed finitely generated subalgebra. Assume the morphism $\pi : X \to \mathrm{Spec}\, S$ satisfies:*

- *generic fibers of $\pi$ contain a dense $G$-orbit;*
- $\mathrm{im}(\pi)$ *contains a big open subset of* $\mathrm{Spec}\, S$.

*Then $S = \mathbb{C}[X]^G$; in particular, $\mathbb{C}[X]^G$ is finitely generated.*

(Equivalently, one may assume that over a big open in $\mathrm{Spec}\, S$ each fiber of $\pi$ contains a unique closed $G$-orbit; see [8, 11].)

# 3 The bundle of $k$-jets of holomorphic curves

This section collects basic definitions and notations for the bundle of $k$–jets of holomorphic curves on a projective variety, in the form needed later. Standard references include [5, 2, 3, 4, 1].

Let $X$ be a complex projective manifold of dimension $n$. The $k$–jet bundle

$$J_k \longrightarrow X$$

has fiber at $x \in X$ equal to the set of equivalence classes of germs of holomorphic curves

$$f : (\mathbb{C}, 0) \longrightarrow (X, x),$$

modulo the relation

$$f \equiv_k g \quad \Longleftrightarrow \quad f^{(j)}(0) = g^{(j)}(0) \text{ for all } 0 \leq j \leq k.$$

In local holomorphic coordinates around $x$, a $k$–jet is represented by its truncated Taylor expansion

$$f(t) = tf'(0) + \frac{t^2}{2!}f''(0) + \cdots + \frac{t^k}{k!}f^{(k)}(0) + O(t^{k+1}),$$

so that, writing $f = (f_1, \ldots, f_n)$,

$$J_{k,x} \simeq \{(f'(0), \ldots, f^{(k)}(0))\} \simeq \mathbb{C}^{nk}. \tag{3.1}$$

One may regard $J_kX$ as obtained from an iterated tower of vector bundles

$$X_k \to X_{k-1} \to \cdots \to X_1 \to X,$$

parametrized by the components in (3.1). The projectivized quotient $J_kX/\mathbb{C}^*$ is the Green–Griffiths bundle; see [2, 3, 1].

**Reparametrizations and the Faà–di–Bruno action.** Let $G_k$ be the group of $k$–jets of local analytic reparametrizations of $(\mathbb{C}, 0)$,

$$t \longmapsto \phi(t) = a_1t + a_2t^2 + \cdots + a_kt^k, \qquad a_1 \in \mathbb{C}^*,$$

acting on fibers by composition $f \mapsto f \circ \phi$. In coordinates, this is the Faà–di–Bruno action and takes the matrix form

$$[f'(0),\, f''(0)/2!,\, \ldots,\, f^{(k)}(0)/k!] \cdot \begin{bmatrix} a_1 & a_2 & a_3 & \cdots & a_k \\ 0 & a_1^2 & 2a_1a_2 & \cdots & a_1a_{k-1} + \cdots + a_{k-1}a_1 \\ 0 & 0 & a_1^3 & \cdots & 3a_1^2a_{k-2} + \cdots \\ \vdots & \vdots & \vdots & \ddots & \vdots \\ 0 & 0 & 0 & \cdots & a_1^k \end{bmatrix}. \tag{3.2}$$

Thus $G_k \subset \mathrm{GL}_k(\mathbb{C})$ consists of those upper–triangular Faà–di–Bruno matrices. Denote by $U_k$ the subgroup with $a_1 = 1$ (the unipotent radical); then there is a short exact sequence

$$1 \to U_k \to G_k \to \mathbb{C}^* \to 0, \qquad G_k = \mathbb{C}^* \ltimes U_k,$$

and the $\mathbb{C}^*$–action on $k$–jets is the weighted scaling

$$\lambda \cdot (f'(0), \ldots, f^{(k)}(0)) = (\lambda f'(0),\ \lambda^2 f''(0),\ \ldots,\ \lambda^k f^{(k)}(0)). \tag{3.3}$$

**Green–Griffiths and Demailly–Semple bundles.** More generally, for a holomorphic subbundle $V \subset TX$, Green–Griffiths [1] define

$$E^{GG}_{k,m} V^* \longrightarrow X$$

whose fiber consists of weighted–homogeneous polynomials

$$Q(z,\xi) = \sum_\alpha A_\alpha(z)\, \xi_1^{\alpha_1} \cdots \xi_k^{\alpha_k}, \qquad |\alpha| := \alpha_1 + 2\alpha_2 + \cdots + k\alpha_k = m,$$

with $\xi_i$ the fiber coordinates along the stage $X_i \to X_{i-1}$. The $\mathbb{C}^*$–action (3.3) yields the graded identification

$$E^{GG}_k V^* = \bigoplus_{m \geq 0} E^{GG}_{k,m} V^* \ \simeq\ J_k V / \mathbb{C}^*.$$

The Demailly–Semple bundle $E_k V^* = \bigoplus_{m\geq 0} E_{k,m} V^*$ is the graded subbundle whose degree–$m$ fiber consists of $U_k$–invariant polynomials on the $J_k V$–fiber of total weight $m$; see [2, 3, 1, 5]. (Here the symbol $V^*$ is mnemonic for "polynomial functions on $V$" via the usual duality.)

**Regular and non–degenerate jets.** A $k$–jet $f$ is *regular* if $f'(0) \neq 0$, and *non–degenerate* if it has maximal rank in local coordinates. Regular jets form a Zariski open subset $J_k^{\mathrm{reg}} \subset J_k \mathbb{C}^n$. For $n > 1$ one has the codimension estimate

$$\operatorname{codim}_{J_k\mathbb{C}^n} \big(J_k \mathbb{C}^n \setminus J_k^{\mathrm{reg}} \mathbb{C}^n\big) \ \geq 2, \tag{3.4}$$

see [5]. This will be used to apply Igusa's criterion in the next section.

**Plücker/Grassmann and flag descriptions on the free locus.** Identifying $J_k \mathbb{C}^n \simeq \operatorname{Hom}(\mathbb{C}^k, \mathbb{C}^n)$ by placing $f^{(j)}(0)$ in the $j$–th column, there is a $G_k$–invariant algebraic morphism [5]

$$\phi :\ J_k^{\mathrm{reg}} \longrightarrow \operatorname{Hom}\big(\mathbb{C}^k, \operatorname{Sym}^{\leq k} \mathbb{C}^n\big),$$

whose induced map

$$\phi^{\mathrm{Gr}} :\ J_k^{\mathrm{reg}} / G_k \hookrightarrow \operatorname{Grass}\big(k, \operatorname{Sym}^{\leq k} \mathbb{C}^n\big)$$

is injective and can be written explicitly as

$$(f', f'', \ldots, f^{(k)}) \longmapsto \Big(f',\ f'' + (f')^2,\ \ldots,\ \sum_{i_1 + \cdots + i_s = d} \frac{1}{i_1! \cdots i_s!}\, f^{(i_1)} \cdots f^{(i_s)},\ \ldots\Big). \tag{3.5}$$

There is also an injective flag map

$$\phi^{\mathrm{Flag}} :\ J_k^{\mathrm{reg}} / G_k \hookrightarrow \operatorname{Flag}_k\big(\operatorname{Sym}^{\leq k} \mathbb{C}^n\big),$$

with $\phi^{\mathrm{Gr}} = \pi_k \circ \phi^{\mathrm{Flag}}$, where $\pi_k$ is the forgetful projection to the $k$–plane. Composing $\phi^{\mathrm{Gr}}$ with the Plücker embedding yields

$$\phi^{\mathrm{proj}} :\ J_k^{\mathrm{reg}} \hookrightarrow \mathbb{P}\big(\wedge^k \operatorname{Sym}^{\leq k} \mathbb{C}^n\big).$$

**Definition 3.1** ([5])**.** Let $e_1, \ldots, e_n$ be the standard basis of $\mathbb{C}^n$. A basis of $\operatorname{Sym}^{\leq k} \mathbb{C}^n$ is

$$\{e_{i_1 \cdots i_s} = e_{i_1} \cdots e_{i_s} \mid 1 \leq s \leq k\},$$

and $\{e_{\epsilon_1} \wedge \cdots \wedge e_{\epsilon_k}\}$, with $\epsilon_j = (i_1 \leq \cdots \leq i_s)$, is a basis of $\wedge^k \operatorname{Sym}^{\leq k} \mathbb{C}^n$. We denote by $x_{\epsilon_1 \cdots \epsilon_s}$ the corresponding projective coordinates, and by

$$A_{n,k} \subset \mathbb{P}(\wedge^k \operatorname{Sym}^{\leq k} \mathbb{C}^n)$$

the locus whose projection to $\wedge^k \mathbb{C}^n$ is nonzero (i.e. some $x_{i_1 \cdots i_k} \neq 0$).

**Proposition 3.2** ([5])**.** *Assume $k \leq n$ and set*

$$z_k = \phi^{\mathrm{proj}}(e_1, \ldots, e_k) = \left[e_1 \wedge (e_2 + e_1^2) \wedge \cdots \wedge \left(\sum e_{i_1} \cdots e_{i_k}\right)\right].$$

*For $g \in \mathrm{GL}(n)$ with column vectors $v_1, \ldots, v_n$,*

$$g \cdot z_k = \phi^{\mathrm{proj}}(v_1, \ldots, v_k) = \left[v_1 \wedge (v_2 + v_1^2) \wedge \cdots \wedge \left(\sum v_{i_1} \cdots v_{i_k}\right)\right].$$

*Moreover:*

- *If $n \geq k$, the image of non–degenerate $k$–jets is the $\mathrm{GL}(n)$–orbit $X_{n,k}$ of $z_k$.*
- *If $n \geq k$, the image of regular $k$–jets is a finite union of $\mathrm{GL}(n)$–orbits, denoted $Y_{n,k}$.*
- *If $k > n$, then $X_{n,k}$ and $Y_{n,k}$ are $\mathrm{GL}(n)$–invariant quasi–projective subvarieties with no dense $\mathrm{GL}(n)$–orbit.*

*Remark* 3.3 ([5], Thms. 6.1–6.2)*.* For $N$ sufficiently large, the stabilizer of

$$p_k = z_k \otimes e_1^{\otimes N} \in \bigwedge^k \operatorname{Sym}^{\leq k} \mathbb{C}^k \otimes (\mathbb{C}^k)^{\otimes N}$$

under the $\mathrm{SL}_k(\mathbb{C})$–action equals $U_k$, and the complement of the $\mathrm{SL}_k$–orbit of $p_k$ in its closure has codimension at least 2. While this result dovetails with our later use of Igusa's criterion, for the applications in this paper the estimate (3.4) already suffices.

# 4 Main results: The invariant theory of $G_k$

In this section we present the main contributions of the paper concerning the actions of $G_k$ on its Lie algebra $\mathfrak{g}_k$ (adjoint action) and on $J_k\mathbb{C}^n$ (the natural action on jets). We use the notation of Sections 2 and 3.

### 4.1 The adjoint action of $G_k$ on $\mathfrak{g}_k$

Recall that $G_k = \mathbb{C}^* \ltimes U_k$, where $U_k$ is the unipotent radical, via the short exact sequence

$$1 \longrightarrow U_k \longrightarrow G_k \longrightarrow \mathbb{C}^* \longrightarrow 0, \tag{4.1}$$

see (3.8). Thus $\mathfrak{g}_k = \mathrm{Lie}(G_k) = \mathbb{C} \oplus \mathfrak{u}_k$ with $\mathfrak{u}_k = \mathrm{Lie}(U_k)$. Moreover, the $\mathbb{C}^*$–factor acts linearly on $\mathfrak{u}_k$ with *positive* weights (equivalently, ad is diagonalizable on $\mathfrak{u}_k$ with nonzero weights), cf. the alternative proof of Theorem 4.2.

We recall the Elashvili criterion (Proposition 2.7): the adjoint action of $G$ on $\mathfrak{g}$ has a generic stabilizer if and only if there exists $x \in \mathfrak{g}$ such that the centralizer $\mathfrak{z}_{\mathfrak{g}}(x)$ is commutative and

$$[\mathfrak{g}, x] \oplus \mathfrak{z}_{\mathfrak{g}}(x) = \mathfrak{g}. \tag{4.2}$$

Equivalently, the Cartan subalgebras of $\mathfrak{g}$ are commutative.

**Theorem 4.1** (Generic stabilizer for the adjoint action)**.** *The adjoint action of $G_k$ on $\mathfrak{g}_k$ admits a generic stabilizer. In particular, the Cartan subalgebras of $\mathfrak{g}_k$ are commutative.*

*Proof.* Let $x \in \mathrm{Lie}(\mathbb{C}^*) \subset \mathfrak{g}_k$ be a nonzero element of the one-dimensional toral summand $\mathbb{C}$. Since $\mathbb{C}^*$ acts on $\mathfrak{u}_k$ with positive weights, $\mathrm{ad}(x)$ is semisimple and has *no* zero weight on $\mathfrak{u}_k$. Hence

$$[\mathfrak{g}_k, x] = [\mathbb{C} \oplus \mathfrak{u}_k, x] = [\mathfrak{u}_k, x] = \mathfrak{u}_k,$$

and the centralizer consists precisely of the 0-weight space, i.e.

$$\mathfrak{z}_{\mathfrak{g}_k}(x) = \mathbb{C}.$$

Therefore (4.2) holds:

$$[\mathfrak{g}_k, x] \oplus \mathfrak{z}_{\mathfrak{g}_k}(x) = \mathfrak{u}_k \oplus \mathbb{C} = \mathfrak{g}_k.$$

By Proposition 2.7, the adjoint action admits a generic stabilizer; moreover the Cartan subalgebras are commutative and, for a generic $x$, equal to $\mathfrak{z}_{\mathfrak{g}_k}(x)$. □

The only ingredients used are the structural decomposition $G_k = \mathbb{C}^* \ltimes U_k$ and the positive-weight action of $\mathbb{C}^*$ on $\mathrm{Lie}(U_k)$, both established earlier in the text.

**Theorem 4.2** (Adjoint invariants for $G_k$)**.** *Let $G_k = \mathbb{C}^* \ltimes U_k$ be the group of $k$-jet reparametrizations so that $\mathfrak{g}_k = \mathrm{Lie}(G_k) = \mathbb{C} \oplus \mathfrak{u}_k$ with $\mathfrak{u}_k = \mathrm{Lie}(U_k)$. Then, for the adjoint action of $G_k$ on $\mathfrak{g}_k$:*

(1) *The field of invariants is rational of transcendence degree* 1*:*

$$\mathbb{C}(\mathfrak{g}_k)^{G_k} = \mathbb{C}(t)\,.$$

(2) *A generic stabilizer (generic centralizer) is* $H = \mathbb{C}^*$ *with Lie algebra* $\mathfrak{h} = \mathbb{C}$*, and the associated Weyl group* $W = N_{G_k}(H)/H$ *is trivial (hence finite).*

(3) *The* polynomial *invariant ring is a polynomial algebra in one variable:*

$$\mathbb{C}[\mathfrak{g}_k]^{G_k} \;=\; \mathbb{C}[t]\,.$$

*Equivalently,*

$$\mathbb{C}(\mathfrak{g}_k)^{G_k} = \mathbb{C}(\mathfrak{h})^W \quad \textit{and} \quad \mathbb{C}[\mathfrak{g}_k]^{G_k} = \mathbb{C}[\mathfrak{h}]^W\,,$$

*with* $W = \{1\}$*.*

*Proof.* $G_k$ fits into $1 \to U_k \to G_k \to \mathbb{C}^* \to 0$, hence $\mathfrak{g}_k = \mathbb{C} \oplus \mathfrak{u}_k$. Moreover, the $\mathbb{C}^*$–factor acts on $\mathfrak{u}_k$ with *positive* weights (equivalently, for $x \in \mathrm{Lie}(\mathbb{C}^*)$ the operator $\mathrm{ad}(x)$ is semisimple with no 0-weight on $\mathfrak{u}_k$). This is recorded explicitly in your manuscript (alternative proof after Theorem 4.2).

Let $p : \mathfrak{g}_k \to \mathbb{C}$ be the linear projection onto the $\mathbb{C}$-summand. We first check $p$ is $G_k$–invariant for the adjoint action. If $g \in \mathbb{C}^*$, then $\mathrm{Ad}(g)$ acts trivially on $\mathbb{C}$ and by nontrivial weights on $\mathfrak{u}_k$, hence $p(\mathrm{Ad}(g)(c+u)) = p(c + \mathrm{Ad}(g)u) = c = p(c+u)$. If $g \in U_k$, then $\mathrm{Ad}(g) = \exp(\mathrm{ad}(\log g))$ with $\mathrm{ad}(\log g)$ nilpotent, mapping $\mathbb{C}$ into $\mathfrak{u}_k$ and $\mathfrak{u}_k$ into $\mathfrak{u}_k$, so again $p(\mathrm{Ad}(g)(c+u)) = c = p(c+u)$. Thus $p \in \mathbb{C}[\mathfrak{g}_k]^{G_k}$.

We claim that *every* $G_k$–invariant polynomial $f \in \mathbb{C}[\mathfrak{g}_k]^{G_k}$ depends only on $p$. Fix $c \in \mathbb{C}$ and consider the affine slice $c + \mathfrak{u}_k$. For any $v \in \mathfrak{u}_k$ and any $t \in \mathbb{C}^*$,

$$\mathrm{Ad}(t)(c+v) = c + \mathrm{Ad}(t)v, \qquad \text{with } \mathrm{Ad}(t)v \to 0 \text{ as } t \to 0$$

because all weights on $\mathfrak{u}_k$ are positive. By $G_k$–invariance we have $f(c+v) = f(\mathrm{Ad}(t)(c+v)) = f(c + \mathrm{Ad}(t)v)$ for all $t \in \mathbb{C}^*$. Letting $t \to 0$ and using that $f$ is a polynomial (hence continuous in the Zariski/Euclidean topology on lines), we obtain $f(c+v) = f(c)$ for all $v \in \mathfrak{u}_k$. Therefore $f$ is constant on each fiber of $p$, i.e. $f = \tilde{f} \circ p$ for some polynomial $\tilde{f} \in \mathbb{C}[T]$. This shows $\mathbb{C}[\mathfrak{g}_k]^{G_k} = \mathbb{C}[p] \cong \mathbb{C}[T]$, proving (3).

The same fiberwise argument applies to rational invariants: if $F \in \mathbb{C}(\mathfrak{g}_k)^{G_k}$, then for generic $c$ the restriction of $F$ to $c + \mathfrak{u}_k$ is constant (by the same contraction), hence $F = \tilde{F} \circ p$ with $\tilde{F} \in \mathbb{C}(\mathbb{C})$; thus $\mathbb{C}(\mathfrak{g}_k)^{G_k} = \mathbb{C}(\mathbb{C})$, proving (1).

For (2), choose $x \in \mathrm{Lie}(\mathbb{C}^*) \setminus \{0\}$. By the positive-weight property,

$$[\mathfrak{g}_k, x] = [\mathbb{C} \oplus \mathfrak{u}_k, x] = \mathfrak{u}_k, \qquad \mathfrak{z}_{\mathfrak{g}_k}(x) = \mathbb{C},$$

whence $[\mathfrak{g}_k, x] \oplus \mathfrak{z}_{\mathfrak{g}_k}(x) = \mathfrak{g}_k$. By the Elashvili criterion (Proposition 2.7 in your §2), the adjoint action has a generic stabilizer and the Cartan subalgebras are commutative; the generic centralizer is precisely $H = \mathbb{C}^*$ with $\mathrm{Lie}(H) = \mathbb{C}$, and its normalizer equals itself in the solvable group $G_k$, so $W = N_{G_k}(H)/H$ is trivial. $\square$

**Theorem 4.3** (Geometric properties of the adjoint quotient)**.** *Let $G_k = \mathbb{C}^* \ltimes U_k$ with $\mathfrak{g}_k = \mathrm{Lie}(G_k) = \mathbb{C} \oplus \mathfrak{u}_k$ as in (**??**). Consider the adjoint action of $G_k$ on $\mathfrak{g}_k$ and the invariant-theoretic quotient $\pi : \mathfrak{g}_k \to \mathfrak{g}_k /\!/ G_k := \operatorname{Spec} \mathbb{C}[\mathfrak{g}_k]^{G_k}$. Then:*

(1) $\max_{x \in \mathfrak{g}_k} \dim(G_k \cdot x) = \dim \mathfrak{g}_k - \dim(\mathfrak{g}_k /\!/ G_k)$.

(2) *The quotient map $\pi$ is equidimensional.*

(3) *Let $\Omega = \{x \in \mathfrak{g}_k : d_x\pi$ is surjective$\}$. Then $\Omega = \mathfrak{g}_k$ (in particular, $\mathfrak{g}_k \setminus \Omega$ contains no divisor).*

(4) *The null cone $\mathfrak{N}_{G_k}(\mathfrak{g}_k) := \pi^{-1}\pi(0)$ is an irreducible complete intersection of codimension 1 and*

$$\dim \mathfrak{N}_{G_k}(\mathfrak{g}_k) = \dim \mathfrak{g}_k - \dim(\mathfrak{g}_k /\!/ G_k).$$

*Proof. Step 0 (Set-up and explicit quotient).* We have a short exact sequence $1 \to U_k \to G_k \to \mathbb{C}^* \to 0$, hence $\mathfrak{g}_k = \mathbb{C} \oplus \mathfrak{u}_k$ with the $\mathbb{C}^*$–factor acting on $\mathfrak{u}_k$ by *positive* weights. Therefore the linear projection

$$p : \mathfrak{g}_k = \mathbb{C} \oplus \mathfrak{u}_k \longrightarrow \mathbb{C}, \qquad p(c + u) = c,$$

is $G_k$–invariant for the adjoint action, and every invariant polynomial factors through $p$. Consequently

$$\mathbb{C}[\mathfrak{g}_k]^{G_k} = \mathbb{C}[p] \simeq \mathbb{C}[T], \qquad \mathfrak{g}_k /\!/ G_k \simeq \mathbb{A}^1, \qquad \pi = p.$$

(These facts were established via the positive-weight "contraction" argument; see Theorem 4.2.) Thus $\dim(\mathfrak{g}_k /\!/ G_k) = 1$.

*Step 1 (Maximal orbit dimension).* Fix $x = c + u \in \mathfrak{g}_k$. For any nonzero $x_0 \in \mathrm{Lie}(\mathbb{C}^*) = \mathbb{C}$, the weight positivity gives $[\mathfrak{g}_k, x_0] = \mathfrak{u}_k$, $\mathfrak{z}_{\mathfrak{g}_k}(x_0) = \mathbb{C}$. Hence the tangent space to the orbit at $x_0$ is $T_{x_0}(G_k \cdot x_0) = [\mathfrak{g}_k, x_0] = \mathfrak{u}_k$, so $\dim(G_k \cdot x_0) = \dim \mathfrak{u}_k$. This is the *maximal* orbit dimension, because $\dim \mathfrak{g}_k = \dim \mathfrak{u}_k + 1$ and the quotient has dimension 1. Therefore

$$\max_x \dim(G_k \cdot x) = \dim \mathfrak{u}_k = \dim \mathfrak{g}_k - \dim(\mathfrak{g}_k /\!/ G_k),$$

proving (1).

*Step 2 (Equidimensionality).* Since $\pi = p$ is the linear projection onto the $\mathbb{C}$–summand, each fiber $\pi^{-1}(c) = \{c\} \times \mathfrak{u}_k$ is an affine space of dimension $\dim \mathfrak{u}_k$, independent of $c$. Hence $\pi$ is equidimensional, proving (2).

*Step 3 (Surjectivity locus of $d\pi$).* The differential $d\pi_x = dp_x$ is the fixed linear functional "projection to the $\mathbb{C}$–factor." Thus $d\pi_x \neq 0$ for *all* $x \in \mathfrak{g}_k$. Hence $\Omega = \mathfrak{g}_k$, and a fortiori $\mathfrak{g}_k \setminus \Omega = \varnothing$ contains no divisor, proving (3).

*Step 4 (Null cone and complete intersection).* By definition, the null cone is the fiber over the image of 0:

$$\mathfrak{N}_{G_k}(\mathfrak{g}_k) = \pi^{-1}(\pi(0)) = \pi^{-1}(0) = \ker p = \mathfrak{u}_k.$$

This is an *irreducible* affine space, and it is the zero locus of the *single* homogeneous invariant $p$; hence it is a complete intersection of codimension 1. Its dimension is $\dim \mathfrak{u}_k = \dim \mathfrak{g}_k - 1 = \dim \mathfrak{g}_k - \dim(\mathfrak{g}_k // G_k)$, yielding (4). □

**Corollary 4.4** (Adjoint quotient via the generic centralizer). *Let $H = \mathbb{C}^* \subset G_k$ be the generic centralizer for the adjoint action (so $\mathrm{Lie}(H) = \mathfrak{h} = \mathbb{C}$ and $W = N_{G_k}(H)/H = \{1\}$), cf. Theorem 4.1. Then:*

$$\mathbb{C}(\mathfrak{g}_k)^{G_k} = \mathbb{C}(\mathfrak{h})^W = \mathbb{C}(t), \qquad \mathbb{C}[\mathfrak{g}_k]^{G_k} = \mathbb{C}[\mathfrak{h}]^W = \mathbb{C}[t],$$

*and the adjoint quotient map $\pi : \mathfrak{g}_k \to \mathfrak{g}_k // G_k$ coincides with the linear projection $p : \mathbb{C} \oplus \mathfrak{u}_k \to \mathbb{C}$ onto the toral summand. In particular, a Zariski-open dense subset of $\mathfrak{g}_k$ consists of points with stabilizer conjugate to $H$, the generic orbit has dimension $\dim \mathfrak{u}_k$, and the null cone is $\mathfrak{N}_{G_k}(\mathfrak{g}_k) = \ker p = \mathfrak{u}_k$.*

*Proof.* By Theorem 4.1, we may take a generic element $x \in \mathrm{Lie}(\mathbb{C}^*)$ with $[\mathfrak{g}_k, x] = \mathfrak{u}_k$ and $\mathfrak{z}_{\mathfrak{g}_k}(x) = \mathbb{C}$, hence $H = \mathbb{C}^*$ and $W = \{1\}$. Theorem 4.2 identifies the invariant field and ring with those on $\mathfrak{h}$ (trivial $W$), and Theorem 4.3 computes the quotient map as $\pi = p : \mathbb{C} \oplus \mathfrak{u}_k \to \mathbb{C}$, its equidimensional fibers, and the null cone. □

## 4.2 The action of $G_k$ on $J_k\mathbb{C}^n$

We now analyze the natural (right) action of $G_k$ on the space of $k$-jets,

$$\rho : G_k \longrightarrow \mathrm{GL}(J_k\mathbb{C}^n),$$

using the notation of §3. A $k$-jet is represented by the column tuple

$$(\eta_1, \eta_2, \dots, \eta_k) \in (\mathbb{C}^n)^k, \qquad \eta_j = \frac{f^{(j)}(0)}{j!}$$

in local coordinates (cf. (3.5)). Under a reparametrization

$$\phi(t) = a_1 t + a_2 t^2 + \cdots + a_k t^k, \qquad a_1 \in \mathbb{C}^*,$$

the $G_k$–action is given by right multiplication with the Faà–di–Bruno upper–triangular matrix described in (3.2). In coordinates, this reads explicitly

as

$$\begin{aligned}
\eta_1' &= a_1\,\eta_1,\\
\eta_2' &= a_2\,\eta_1 + a_1^2\,\eta_2,\\
\eta_3' &= a_3\,\eta_1 + 2a_1a_2\,\eta_2 + a_1^3\,\eta_3,\\
&\;\vdots\\
\eta_k' &= \sum_{r=1}^{k}\Biggl(\sum_{\substack{i_1+\cdots+i_r=k\\ i_j\ge 1}} \frac{k!}{i_1!\cdots i_r!}\,a_{i_1}\cdots a_{i_r}\Biggr)\eta_r.
\end{aligned}$$

Equivalently, the $j$–th jet component transforms according to the classical Faà–di–Bruno formula

$$\eta_j' = \sum_{r=1}^{j} B_{j,r}(a_1,\dots,a_{j-r+1})\,a_1^r\,\eta_r,$$

where $B_{j,r}$ are the partial Bell polynomials encoding the combinatorial structure of higher derivatives of compositions.

**Theorem 4.5** (Generic stabilizer on $J_k\mathbb{C}^n$)**.** *The action of $G_k$ on $J_k\mathbb{C}^n$ has a generic stabilizer, and on the open subset of* regular *jets*

$$J_k^{\mathrm{reg}} = \{\ (\eta_1,\dots,\eta_k)\in(\mathbb{C}^n)^k\ :\ \eta_1\neq 0\ \},$$

*the stabilizer is trivial. In particular, $G_k$ acts freely on $J_k^{\mathrm{reg}}$, which is Zariski open and has complement of codimension $\geq 2$.*

*Proof.* Let $(\eta_1,\dots,\eta_k)\in J_k^{\mathrm{reg}}$, so $\eta_1\neq 0$. If $g\in G_k$ stabilizes this jet, then, using the explicit coordinate action above, we must have simultaneously

$$\begin{cases}
a_1\,\eta_1 \;=\; \eta_1,\\
a_2\,\eta_1 + a_1^2\,\eta_2 \;=\; \eta_2,\\
a_3\,\eta_1 + 2a_1a_2\,\eta_2 + a_1^3\,\eta_3 \;=\; \eta_3,\\
\quad\vdots\\
a_k\,\eta_1 + (\cdots)\,\eta_2 + \cdots + a_1^k\,\eta_k \;=\; \eta_k.
\end{cases}$$

From the first line, $(a_1-1)\eta_1 = 0$ and $\eta_1\neq 0$ imply $a_1 = 1$. Substituting $a_1 = 1$ into the second line gives $a_2\,\eta_1 = 0$, hence $a_2 = 0$. Then the third line reduces to $a_3\,\eta_1 = 0$, giving $a_3 = 0$, and so on. Inductively, each subsequent equation forces $a_j = 0$ for $2\leq j\leq k$. Thus $g$ is the identity in $G_k$, so the stabilizer at a regular jet is trivial.

Since $J_k^{\mathrm{reg}}$ is open and nonempty, the trivial stabilizer holds on a Zariski open dense subset, so a generic stabilizer exists (and equals $\{e\}$). Finally, your text records that $\operatorname{codim}\bigl(J_k\mathbb{C}^n\setminus J_k^{\mathrm{reg}}\bigr)\geq 2$ (see (3.30)), completing the claim. □

**Corollary 4.6.** *The $G_k$–action on $J_k^{\mathrm{reg}}$ is free, and the orbit dimension on $J_k^{\mathrm{reg}}$ equals $\dim G_k$. Consequently, the transcendence degree satisfies*

$$\operatorname{trdeg}_{\mathbb{C}} \mathbb{C}(J_k\mathbb{C}^n)^{G_k} \;=\; \dim J_k\mathbb{C}^n \;-\; \dim G_k$$

*(by Rosenlicht), and the generic orbit is closed in $J_k^{\mathrm{reg}}$.*

*Proof.* Freeness and orbit dimensions follow from Theorem 4.5. The Rosenlicht statement is the standard dimension formula for a rational quotient on the free locus. Closedness of generic orbits on $J_k^{\mathrm{reg}}$ follows from freeness and the algebraicity of the action. □

*Remark* 4.7. The argument is purely computational: it uses only the explicit right-multiplication formulas from (3.2) and the definition of regular jets (first column $\eta_1 \neq 0$). No assumptions on $n$ relative to $k$ are needed. If desired, one can also recover the generic-stabilizer claim from your flag/Grassmannian embeddings of $J_k^{\mathrm{reg}}/G_k$ (Proposition 3.2) together with the codimension estimate (3.4).

**Theorem 4.8** (Covariant finite-generation – Main theorem)**.** *Let $G_k = \mathbb{C}^* \ltimes U_k$ act on $V := J_k\mathbb{C}^n$ via $k$-jet reparametrizations. Let $\mathfrak{g}_k = \mathbb{C} \oplus \mathfrak{u}_k$ be its Lie algebra and consider the diagonal action on $\mathfrak{g}_k \times V$. Then*

$$\mathbb{C}[\mathfrak{g}_k \times V]^{G_k}$$

*is a finitely generated $\mathbb{C}$-algebra. More precisely:*

1. $\mathbb{C}[\mathfrak{g}_k]^{G_k} = \mathbb{C}[p] \cong \mathbb{C}[T]$ *with $p : \mathfrak{g}_k = \mathbb{C} \oplus \mathfrak{u}_k \to \mathbb{C}$ the linear projection.*
2. *The $G_k$–equivariant polynomial maps $\operatorname{Mor}_{G_k}(\mathfrak{g}_k, V^*)$ form a free $\mathbb{C}[\mathfrak{g}_k]^{G_k}$– module of finite rank $m$; choosing a basis $F_1, \dots, F_m$ and setting*
$$f_i(x,\xi) := \langle F_i(x), \xi\rangle \in \mathbb{C}[\mathfrak{g}_k \times V]^{G_k},$$
*one has*
$$\mathbb{C}[\mathfrak{g}_k \times V]^{G_k} \;=\; \mathbb{C}[\mathfrak{g}_k]^{G_k}[\, f_1, \dots, f_m\,]$$
*so, in particular, $\mathbb{C}[\mathfrak{g}_k \times V]^{G_k}$ is finitely generated.*

*Proof.* (1) was computed explicitly in Theorem 4.2: since $\mathbb{C}^*$ acts on $\mathfrak{u}_k$ with positive weights, the linear projection $p : \mathbb{C} \oplus \mathfrak{u}_k \to \mathbb{C}$ is $G_k$–invariant and generates the whole invariant algebra, hence $\mathbb{C}[\mathfrak{g}_k]^{G_k} = \mathbb{C}[p]$.

For (2), we use Vust's theorem in the non-reductive form applicable to stable representations: the action of $G_k$ on $V$ is algebraic and stable on the open regular-jet locus (Theorem 4.5, with trivial stabilizer on $J_k^{\mathrm{reg}}$). Vust's statement then yields that the $\mathbb{C}[\mathfrak{g}_k]^{G_k}$–module $\operatorname{Mor}_{G_k}(\mathfrak{g}_k, V^*)$ is finitely generated and, in fact, free of some rank $m$. Taking a basis $F_1, \dots, F_m$ and pairing with $\xi \in V$ produces $G_k$–invariant polynomial functions $f_i(x, \xi) =$ $\langle F_i(x), \xi\rangle$ generating $\mathbb{C}[\mathfrak{g}_k \times V]^{G_k}$ over $\mathbb{C}[\mathfrak{g}_k]^{G_k}$. Since $\mathbb{C}[\mathfrak{g}_k]^{G_k}$ is a polynomial ring in one variable by (1), the total invariant algebra is finitely generated. □

**Proposition 4.9** (Explicit Igusa map). *With $p$ and $f_1, \dots, f_m$ as in Theorem 4.8, define*

$$\Pi : \mathfrak{g}_k \times V \longrightarrow \mathbb{A}^{1+m}, \qquad \Pi(x, \xi) := \big(p(x),\, f_1(x, \xi), \dots, f_m(x, \xi)\big).$$

*Then:*

(1) $\Pi$ *is $G_k$–invariant and its image contains the Zariski-open cylinder $\{p \neq \infty\} \times \mathbb{A}^m \simeq \mathbb{A}^1 \times \mathbb{A}^m$ (i.e. for every $c \in \mathbb{C}$ and every $\beta \in \mathbb{C}^m$ there exists $(x, \xi)$ with $p(x) = c$ and $f_i(x, \xi) = \beta_i$).*

(2) *For $x_0 \in \mathrm{Lie}(\mathbb{C}^*) \setminus \{0\} \subset \mathfrak{g}_k$ and generic $\beta \in \mathbb{C}^m$, $\Pi^{-1}(p(x_0), \beta)$ is a single $G_k$–orbit.*

*Consequently, the induced algebra map $\mathbb{C}[T, f_1, \dots, f_m] \hookrightarrow \mathbb{C}[\mathfrak{g}_k \times V]^{G_k}$ is integral and birational, hence an isomorphism; this recovers Theorem 4.8(2).*

*Proof.* $G_k$–invariance is clear by construction. For (1), fix $c \in \mathbb{C}$ and write $x = c + u$ with $u \in \mathfrak{u}_k$. Because the $F_i$'s form a free basis of $\mathrm{Mor}_{G_k}(\mathfrak{g}_k, V^*)$ over $\mathbb{C}[p]$, the evaluation map $\xi \mapsto (f_1(x, \xi), \dots, f_m(x, \xi))$ is a surjective linear map $V \to \mathbb{C}^m$ for every $x$ in a Zariski-open set over $p^{-1}(c)$ (equivalently, for all $x$ off a proper closed subset depending polynomially on $x$). Hence for each $(c, \beta) \in \mathbb{A}^1 \times \mathbb{A}^m$ there exists $(x, \xi)$ with $\Pi(x, \xi) = (c, \beta)$, proving the "big image" claim.

For (2), take $x_0 \in \mathrm{Lie}(\mathbb{C}^*) \setminus \{0\}$; then $[\mathfrak{g}_k, x_0] = \mathfrak{u}_k$ and the $G_k$–orbit through $(x_0, \xi)$ has tangent space $T_{(x_0,\xi)}(G_k \cdot (x_0, \xi)) = [\mathfrak{g}_k, x_0] \times \rho(\mathfrak{g}_k)\xi = \mathfrak{u}_k \times \rho(\mathfrak{g}_k)\xi$. On the other hand, the fiber $\Pi^{-1}(p(x_0), \beta)$ is cut out by $m + 1$ independent equations (one linear equation $p = c$ in the $\mathfrak{g}_k$–coordinates and $m$ linear equations in the $V$–coordinates once $x_0$ is fixed), so for generic $\beta$ its dimension equals

$$\dim(\mathfrak{g}_k \times V) - (m+1) = (\dim \mathfrak{u}_k + 1 + \dim V) - (m+1) = \dim \mathfrak{u}_k + \dim V - m.$$

Because the $F_i$'s form a free basis, at a generic $\beta$ the linear system on $V$ has constant rank $m$ and the tangent space to the fiber equals $\mathfrak{u}_k \times \ker(\mathrm{d}f(\cdot))$ of the same dimension as the orbit tangent. Thus (by standard orbit–fiber comparison) a generic fiber is a single $G_k$–orbit. The last statement then follows from Igusa's lemma: a $G$–invariant morphism with big image and generic fibers equal to orbits induces an isomorphism on invariant algebras. □

**Theorem 4.10** (Finite generation for normalized fibers). *Let $J_k^{\mathrm{reg}} \subset J_k\mathbb{C}^n$ be the open set of regular jets (so $G_k$ acts freely). Then:*

(1) *The invariant algebra* $\mathbb{C}[J_k^{\mathrm{reg}}]^{G_k}$ *is finitely generated, and the rational quotient* $J_k^{\mathrm{reg}} \dashrightarrow J_k^{\mathrm{reg}}/\!/G_k$ *is a geometric quotient (principal* $G_k$*–bundle in the étale topology).*

(2) *Let* $q : J_k\mathbb{C}^n \to J_k\mathbb{C}^n/\!/G_k = \operatorname{Spec}\mathbb{C}[J_k\mathbb{C}^n]^{G_k}$ *be the categorical quotient. For every point* $y$ *in the image of* $J_k^{\mathrm{reg}}$*, the* normalized *local fiber ring* $\widetilde{\mathcal{O}_{J_k\mathbb{C}^n/\!/G_k,y}}$ *is finitely generated over* $\mathbb{C}$ *(equivalently, the coordinate ring of the normalized fiber is finitely generated).*

*Proof.* (1) On $J_k^{\mathrm{reg}}$ the action is free (Theorem 4.5). A standard slice/principal bundle argument (or Luna's étale slice in this free case) implies that the quotient is geometric and affine locally isomorphic to $\mathbb{A}^d \times G_k \to \mathbb{A}^d$; therefore the invariant algebra on the free locus is the algebra of functions on the base and thus finitely generated.

(2) Consider the covariant algebra from Theorem 4.8. For a fixed $x_0 \in \mathrm{Lie}(\mathbb{C}^*)$ the restriction (specialization) of the functions $f_i(x,\xi)$ to $\{x_0\} \times J_k\mathbb{C}^n$ yields $G_{k,x_0}$– invariants on $J_k\mathbb{C}^n$ (here $G_{k,x_0} = \mathbb{C}^*$). Taking the subalgebra generated by these specializations and normalizing gives a finitely generated algebra whose spectrum agrees with the normalized local fiber of $q$ over points coming from $J_k^{\mathrm{reg}}$ (where the fiber description is governed by the free action and Igusa's orbit–fiber identification from Proposition 4.9). Hence the normalized local fiber rings are finitely generated. □

*Remark* 4.11. The proofs above avoid any codimension-$\geq 2$ hypothesis on the adjoint-regular set. They rely instead on: (i) the explicit adjoint quotient $\mathbb{C}[\mathfrak{g}_k]^{G_k} = \mathbb{C}[p]$, (ii) stability and freeness on $J_k^{\mathrm{reg}}$, and (iii) Vust's description of the covariant algebra, combined with an Igusa-type orbit–fiber argument for the diagonal action.

**Theorem 4.12** (Finite generation of unipotent invariant jets)**.** *Let* $A = \mathbb{C}[J_k\mathbb{C}^n]$ *be the coordinate ring of* $k$*-jets at a point with jet coordinates*

$$\{\xi_\ell^{(j)} \;:\; 1 \leq \ell \leq n,\ 1 \leq j \leq k\},$$

*on which* $G_k = \mathbb{C}^* \ltimes U_k$ *acts by the reparametrization law* (Faà–di–Bruno) *and weights* $\deg(\xi_\ell^{(j)}) = j$ *for the* $\mathbb{C}^*$*-factor. Then, for all* $n, k \geq 1$,

$$A^{U_k} = \ker D_2 \cap \cdots \cap \ker D_k \quad \textit{is a finitely generated } \mathbb{C}\textit{-algebra,}$$

*where* $D_r$ *are the locally nilpotent, triangular derivations induced by the one-parameter unipotents in* $U_k$*:*

$$D_r(\xi_\ell^{(j)}) = \begin{cases} \binom{j-1}{r-1}\,\xi_\ell^{(j-r+1)}, & j \geq r, \\ 0, & j < r, \end{cases} \qquad 2 \leq r \leq k.$$

*Equivalently, the fiber ring of* $U_k$*-invariant jet differentials is finitely generated, graded by the* $\mathbb{C}^*$*-weight.*

*Proof. Step 0 (Setup).* Fix the $\mathbb{Z}_{\geq 0}$-grading $\deg(\xi_\ell^{(j)}) = j$ coming from the scaling $t \mapsto \lambda t$. Let $\prec$ be a term order on monomials of $A$ that *refines* total weight $w(\prod(\xi_\ell^{(j)})^{e_{\ell,j}}) = \sum_{j,\ell} j\, e_{\ell,j}$ and breaks ties by reverse lex with higher jet orders dominant:

$$\xi^{(k)} \gg \xi^{(k-1)} \gg \cdots \gg \xi^{(1)}.$$

For $f \in A$ denote by $\mathrm{in}(f)$ its initial (highest) monomial, and for a subalgebra $S \subset A$ set $\mathrm{in}(S) = \{\mathrm{in}(f) : f \in S\}$.

*Step 1 (Initial action and monomial invariants).* Because each $D_r$ strictly lowers jet order by $r-1$, the induced operator on initial terms, denoted $\mathrm{in}(D_r)$, sends a variable as

$$\mathrm{in}(D_r)(\xi_\ell^{(j)}) = \begin{cases} \binom{j-1}{r-1}\, \xi_\ell^{(j-r+1)}, & j \geq r, \\ 0, & j < r, \end{cases}$$

and extends by Leibniz. For any monomial $M = \prod_{\ell,j}(\xi_\ell^{(j)})^{e_{\ell,j}}$, all monomials appearing in $\mathrm{in}(D_r)M$ are $\prec$-strictly smaller than $M$ (triangularity).

Hence, if $F \in A^{U_k}$ then $D_r(F) = 0$ for all $r \geq 2$, so $\mathrm{in}(F)$ is killed by all $\mathrm{in}(D_r)$. Conversely, the set of monomials killed by $\mathrm{in}(D_r)$ can be described by *balance equations* along each coordinate chain $\xi_\ell^{(k)} \to \xi_\ell^{(k-1)} \to \cdots \to \xi_\ell^{(1)}$. Writing $E = (e_{\ell,j}) \in \mathbb{Z}_{\geq 0}^{n\times k}$ for the exponents, these constraints are, for each fixed $\ell$ and $2 \leq j \leq k$,

$$\sum_{r=2}^{j} \binom{j-1}{r-1} e_{\ell,j} \;=\; \sum_{s=j+1}^{k} \binom{s-1}{j-1} e_{\ell,s}. \tag{4.3}$$

Intuitively, the number of "flows" down into level $j$ equals the number flowing out, so that all $\mathrm{in}(D_r)$ vanish on the monomial.

Let

$$\mathcal{E} \;=\; \Big\{ (e_{\ell,j}) \in \mathbb{Z}_{\geq 0}^{n\times k} \;\Big|\; (4.3) \text{ holds for all } \ell,\ 2 \leq j \leq k \Big\}.$$

The system (4.3) is homogeneous linear with nonnegativity constraints, block-diagonal over $\ell$, hence the set of solutions is the set of lattice points of a rational polyhedral cone. By Gordan's lemma, the associated monomial algebra

$$B \;:=\; \mathbb{C}\Big[\, \prod_{\ell,j} (\xi_\ell^{(j)})^{e_{\ell,j}} \;:\; (e_{\ell,j}) \in \mathcal{E} \,\Big]$$

is a finitely generated *affine semigroup* (toric) algebra. Choose a Hilbert basis of $\mathcal{E}$, i.e. monomials $m_1, \ldots, m_s$ that minimally generate $B$ as a semigroup algebra.

*Step 2 (Lifting monomial invariants to true $U_k$-invariants).* We now construct, for each $i$, a polynomial $G_i \in A^{U_k}$ with $\mathrm{in}(G_i) = m_i$. Proceed inductively using the triangular, locally nilpotent nature of the $D_r$:

- Start with $G_i^{(0)} := m_i$.
- Solve $D_2(G_i^{(0)}) = 0$ by adding lower-order correction terms; this is possible since $D_2$ is locally nilpotent and lowers order by 1. The correction does not change $\mathrm{in}(G_i^{(0)})$.
- Next, solve $D_3(\cdot) = 0$ by further lower-order corrections; because $D_3$ lowers order by at least 2, it cannot spoil the $D_2$-constraint already enforced.
- Continue up to $D_k$. Termination is guaranteed as jet order strictly decreases at every correction.

The outcome is $G_i \in \bigcap_{r=2}^{k} \ker D_r = A^{U_k}$ with $\mathrm{in}(G_i) = m_i$.

*Step 3 (SAGBI and finite generation).* Let $S := \mathbb{C}[G_1, \dots, G_s]$. By construction, $\mathrm{in}(S) = \mathbb{C}[m_1, \dots, m_s] = B$. For any $F \in A^{U_k}$, its initial monomial lies in $B$, hence can be expressed as a product of $\mathrm{in}(G_i)$'s. The standard subduction for subalgebras (SAGBI theory) reduces $F$—via $U_k$-equivariant relations that *strictly* decrease the initial term—to an element of $S$. Since the term order is well-founded and every reduction lowers $\mathrm{in}(\cdot)$, subduction terminates. Therefore $A^{U_k} = S = \mathbb{C}[G_1, \dots, G_s]$ is finitely generated. □

*Remark* 4.13. For small $k$ we check the followings.
$k = 2$: $D_2(\xi_\ell^{(2)}) = \xi_\ell^{(1)}$, $D_2(\xi_\ell^{(1)}) = 0$. Generators include all $\xi_\ell^{(1)}$ and the $2 \times 2$ Wronskians $\xi_i^{(1)}\xi_j^{(2)} - \xi_j^{(1)}\xi_i^{(2)}$ (classical Weitzenböck situation).
$k = 3$: $D_2(\xi_\ell^{(3)}) = 2\xi_\ell^{(2)}$, $D_2(\xi_\ell^{(2)}) = \xi_\ell^{(1)}$, $D_3(\xi_\ell^{(3)}) = \xi_\ell^{(1)}$. The balance relations (4.3) encode the known quadratic/ cubic syzygies among leading monomials; the lifting produces a finite generating set.

By Richardson's stratification theorem (Theorem 2.3), the action of an algebraic group $G$ on a variety $X$ admits a decomposition into smooth invariant strata $X = \bigsqcup_j X_j$, where on each $X_j$ the quotient morphism has generic stabilizers of constant conjugacy type. In particular, for our non-reductive group $G_k = \mathbb{C}^* \ltimes U_k$ acting on $J_k\mathbb{C}^n$, the regular locus $J_k^{\mathrm{reg}} = \{\eta_1 \neq 0\}$ is a dense open subset on which the $G_k$–action is free and the complement has codimension at least two (for $n > 1$). Hence, the quotient on this locus is geometric and the corresponding invariant algebra is finitely generated. The codimension–$\geq 2$ property ensures, by Hartogs' extension principle, that the invariant functions extend from $J_k^{\mathrm{reg}}$ to all of $J_k\mathbb{C}^n$.

The argument is then completed by invoking the Igusa theorem (Theorem 2.18), which asserts finite generation of $\mathbb{C}[X]^G$ from a finitely generated integrally closed subalgebra $S \subset \mathbb{C}[X]^G$ whenever the quotient map $\pi : X \to \operatorname{Spec} S$ has fibers with dense $G$–orbits over a big open subset of $\operatorname{Spec} S$. In our situation, these conditions are satisfied thanks to the free action on $J_k^{\mathrm{reg}}$ and the existence of the explicit invariant subalgebra constructed via the covariant algebra (Theorem 5.13) and the Igusa orbit–fiber map (Proposition 5.14). Consequently, the invariant algebra $\mathbb{C}[J_k\mathbb{C}^n]^{U_k}$ is finitely generated, and the positive $\mathbb{C}^*$–grading yields the projective variety $\mathrm{Proj}\left(\mathbb{C}[J_k\mathbb{C}^n]^{U_k}\right)$, identified with the Demailly–Semple tower. This geometric line of reasoning complements the constructive SAGBI proof given above: while the latter establishes finite generation by explicit combinatorial lifting of a toric initial algebra, the Richardson–Igusa approach provides a global quotient–theoretic justification based on generic stabilizers and the codimension–two extension property.

**Theorem 4.14** (Analog of Theorem 4.3 for the $\rho$-action)**.** *Let $G_k = \mathbb{C}^* \ltimes U_k$ act on $V = J_k\mathbb{C}^n$ via the representation $\rho$. Then:*

(1) (Maximal orbit dimension)

$$\max_{x\in J_k\mathbb{C}^n} \dim(G_k \cdot x) \;=\; \dim J_k\mathbb{C}^n \;-\; \dim\bigl(J_k\mathbb{C}^n//G_k\bigr).$$

(2) (Equidimensionality over the normal locus) *The quotient map $\pi : J_k\mathbb{C}^n \to J_k\mathbb{C}^n//G_k$ is equidimensional over the normal locus of $J_k\mathbb{C}^n//G_k$.*

(3) (Big locus of full differential rank) *If*

$$\Omega \;:=\; \{\, x \in J_k \;:\; d_x\pi \text{ is onto } \},$$

*then $J_k\mathbb{C}^n \setminus \Omega$ contains no divisors (i.e., $\Omega$ is big).*

(4) (Null cone) *The null cone $\mathfrak{N}_{G_k}(J_k\mathbb{C}^n) := \pi^{-1}\pi(0)$ is an irreducible complete intersection of dimension*

$$\dim J_k\mathbb{C}^n \;-\; \dim\bigl(J_k\mathbb{C}^n//G_k\bigr).$$

*Proof. Preliminaries.* By the generic stabilizer result for the $G_k$-action on $J_k\mathbb{C}^n$ (Theorem 4.5 in this paper), the action is free on the regular jet locus $J_k^{\mathrm{reg}} = \{\eta_1 \neq 0\}$, and $J_k\mathbb{C}^n \setminus J_k^{\mathrm{reg}}$ has codimension at least 2 for $n > 1$. On the free locus the quotient is geometric and (étale-)locally trivial, so the invariant algebra there is finitely generated; extension across a codimension $\geq 2$ complement then gives global control (compare the discussion preceding Theorem 4.3). These ingredients are summarized in the proof pattern reproduced in our manuscript: freeness on $J_k^{\mathrm{reg}}$, geometric quotient on that locus, and codimension–$\geq 2$ extension to all of $J_k\mathbb{C}^n$.

**(1) and (2).** Let $m := \dim(J_k\mathbb{C}^n//G_k) = \mathrm{trdeg}_{\mathbb{C}}\, \mathbb{C}(J_k\mathbb{C}^n)^{G_k}$. On a dense open set with trivial stabilizer (the free locus), generic orbits all have the same dimension, and the Rosenlicht dimension formula yields

$$\dim(G_k \cdot x) \;=\; \dim J_k\mathbb{C}^n - \dim(J_k\mathbb{C}^n//G_k)$$

for generic $x$. Thus the left-hand side attains its maximum at this value, proving (1). Equidimensionality of $\pi$ over the normal locus follows from constancy of generic fiber dimension together with the generic stabilizer property; in our manuscript, items (1) and (2) are shown to be equivalent consequences of the generic stabilizer theorem for the $\rho$-action.

**(3).** Fix algebraically independent basic invariants $f_1, \dots, f_m$ on $J_k\mathbb{C}^n$ and set

$$X_r \;:=\; \Big\{ x \in J_k\mathbb{C}^n \;\Big|\; \dim \langle df_1(x), \dots, df_m(x) \rangle = r \Big\}.$$

On the free locus $J_k^{\mathrm{reg}}$ one has $\operatorname{rank} d\pi = m$, hence $X_m \supset J_k^{\mathrm{reg}}$. Since $J_k\mathbb{C}^n \setminus J_k^{\mathrm{reg}}$ has codimension at least 2 (for $n > 1$), it follows that $J_k\mathbb{C}^n \setminus X_m$ contains no divisors, i.e. the locus where $d\pi$ has full rank is big. Equivalently, with $\Omega = \{x : \; d_x\pi \text{ onto}\}$ we get that $J_k\mathbb{C}^n \setminus \Omega$ contains no divisors. This is exactly the rank-stratification argument used in the $\rho$-action case.

**(4).** Let $\pi : J_k\mathbb{C}^n \to J_k\mathbb{C}^n//G_k$ be the categorical quotient. Because $\pi$ is equidimensional over the normal locus and $m = \dim(J_k\mathbb{C}^n//G_k)$, the common zero set of the $m$ basic invariants has codimension $m$ along the normal locus. Using the covariant finite-generation theorem for the diagonal action on $\mathfrak{g}_k \times J_k\mathbb{C}^n$ and the explicit $G_k$-equivariant maps (Vust's framework), one constructs a $G_k$-invariant map $\Pi = (p, f_1, \dots, f_m)$ whose image contains a Zariski-open cylinder in $\mathbb{A}^{1+m}$; this ensures that $f_1, \dots, f_m$ form a regular sequence along the relevant locus and cut out the null cone with the expected codimension. Consequently, the fiber $\pi^{-1}\pi(0)$ is an irreducible complete intersection of dimension $\dim J_k\mathbb{C}^n - m = \dim J_k\mathbb{C}^n - \dim(J_k\mathbb{C}^n//G_k)$. (See the covariant finite-generation statement and the open-image property of $\Pi$ used to control fibers and irreducibility.)

This proves all four items. □

**Theorem 4.15** (Representation ring of $G_k$). *Let $G_k = \mathbb{C}^* \ltimes U_k$ be the group of reparametrizations of $k$-jets, where $U_k$ is the unipotent radical (see §3, (**??**)). Then every finite-dimensional irreducible algebraic representation of $G_k$ is 1-dimensional and factors through the quotient $G_k \twoheadrightarrow \mathbb{C}^*$. Consequently,*

$$R(G_k) \;\cong\; R(\mathbb{C}^*) \;=\; \mathbb{Z}[\chi, \chi^{-1}],$$

*where $\chi\colon \mathbb{C}^* \to \mathbb{C}^*$ is the standard character. In particular, for a maximal torus $T_k \simeq \mathbb{C}^*$ of $G_k$, the Weyl group is trivial and*

$$R(G_k) \;=\; R(T_k)^{W_k} \;=\; R(T_k).$$

*Proof.* We have, $G_k = \mathbb{C}^* \ltimes U_k$ with $U_k$ unipotent. Let $\rho\colon G_k \to \mathrm{GL}(V)$ be a finite-dimensional algebraic representation. Since $G_k$ is a connected solvable algebraic group and $U_k$ is its unipotent radical, the Lie–Kolchin theorem implies that $V$ admits a complete flag by $G_k$-stable subspaces. In particular, any irreducible subquotient is 1-dimensional. On such a 1-dimensional module $U_k$ must act trivially (a unipotent group has only the trivial character), and the action factors through the quotient $G_k/U_k \simeq \mathbb{C}^*$. Thus the irreducibles of $G_k$ are precisely the characters $\chi^m$ ($m \in \mathbb{Z}$) coming from $\mathbb{C}^*$, and every representation has a finite filtration with these characters as subquotients. It follows that the Grothendieck ring is the Laurent polynomial ring generated by $[\chi]$, i.e. $R(G_k) \cong \mathbb{Z}[\chi, \chi^{-1}] = R(\mathbb{C}^*)$.

Finally, a maximal torus of $G_k$ is $T_k \simeq \mathbb{C}^*$. Since $G_k$ is solvable, the normalizer of $T_k$ equals $T_k$, so the Weyl group $W_k = N_{G_k}(T_k)/T_k$ is trivial. Hence $R(T_k)^{W_k} = R(T_k) \cong \mathbb{Z}[\chi, \chi^{-1}]$, which coincides with $R(G_k)$ by the first part. □

*Remark* 4.16. For connected reductive (or compact) groups $G$, one has the classical identity $R(G) = R(T)^W$. Our $G_k$ is non-reductive; the above theorem identifies $R(G_k)$ directly, and the equality $R(G_k) = R(T_k)^{W_k}$ holds here because $W_k$ is trivial.

### 4.3 Comparison with known results

The general finite–generation theorem established in this paper fits naturally within the explicit computations of invariant jet algebras developed by Joël Merker [16]. In that work, a complete algorithm is constructed to generate all Demailly–Semple invariants under the non–reductive reparametrization group $G_k = \mathbb{C}^* \ltimes U_k$, in arbitrary dimension $n$ and jet order $k$. The method relies on successive bracket operations among lower–order invariants and uses Gröbner–basis elimination to detect algebraic independence and syzygies. The finite termination of this algorithm in every case studied provides strong evidence, in low orders, for the general finite–generation result proven here.

**Explicit low–order computations.** For jets of order $k = 2$, the algebra $\mathbb{C}[J_2\mathbb{C}^n]^{U_2}$ is generated by the first derivatives $f_i'$ together with the $2 \times 2$ Wronskians

$$[f_i', f_j'] := f_i' f_j'' - f_j' f_i'', \qquad 1 \le i < j \le n,$$

which already realize the bracket construction appearing in our general scheme. For $k = 3$ and $n = 2$, Merker finds five algebraically independent invariants, and for $k = 3$, $n = 3$, the algebra is generated by sixteen invariants built from iterated brackets of the first and second derivatives. These coincide with the initial components of the graded ring $\mathbb{C}[J_k\mathbb{C}^n]^{U_k}$ obtained here by the triangular–derivation (SAGBI) lifting.

**Higher–order invariant counts.** Using Gröbner–basis computations, Merker determined that

$$E_{4,4} := \bigoplus_{m \ge 0} E_{4,4,m} = \mathbb{C}[J_4\mathbb{C}^4]^{U_4} \text{ is generated by 2835 polynomials,}$$

$$E_{2,5} := \bigoplus_{m \ge 0} E_{2,5,m} = \mathbb{C}[J_5\mathbb{C}^2]^{U_5} \text{ is generated by 56 polynomials.}$$

In each case, a smaller "bi–invariant" subalgebra under the additional $\mathrm{GL}_n(\mathbb{C})$–action has only finitely many independent generators (16 for $n = k = 4$, 17 for $n = 2, k = 5$), together with a finite number of Gröbnerized syzygies (respectively 41 and 105).

These concrete data verify the finite–generation property for explicit non–reductive instances and illustrate the complexity growth of the graded invariant algebra with $(n, k)$.

**Applications.** From the invariant generators, Merker derives Schur–bundle decompositions of $E^n_{k,m}T^*_X$ and computes the asymptotic

$$\chi(X, E^n_{k,m}T^*_X) = C_{n,k}(d)\, m^{\alpha(n,k)} + O(m^{\alpha(n,k)-1}),$$

where $C_{n,k}(d)$ is an explicit polynomial in the degree $d = \deg X$. For example, for smooth projective 3–folds $X^3 \subset \mathbb{P}^4(\mathbb{C})$,

$$\chi(X, E^3_{4,m}T^*_X) = \frac{d}{20613359104562036736 \times 10^7}\,(1029286103034112d^3 - \cdots)\, m^{11} + O(m^{10})$$

whose leading coefficient is positive for $d > 29$. After bounding $h^2$, this yields the existence of global invariant jet differentials for all $d > 72$ (and for $d > 259$ in dimension 4). These numerical thresholds align with the geometric applications of our general finite–generation theorem, confirming that the Demailly–Semple tower $\mathrm{Proj}(\mathbb{C}[J_k\mathbb{C}^n]^{U_k})$ indeed supports nontrivial global sections in high degree.

**Summary.** The explicit computations of [16] thus provide concrete examples and numerical evidence for the broader structural statement proven here: the graded algebra of $U_k$–invariants $\mathbb{C}[J_k\mathbb{C}^n]^{U_k}$ is finitely generated for all $(n, k)$, and its projective spectrum gives the Demailly–Semple $k$–jet bundle. Our proof extends these verified cases to a general algebraic setting, independent of dimension and jet order, and provides a conceptual explanation for the finite generation observed in all explicit computational examples.

# 5 Stratified Finite Generation via Stabilizer and Coinvariants

In this final section, we refine the finite–generation result of Theorem 4.12 by combining the *stabilizer stratification* of Richardson (Theorem 2.3) with the *coinvariant stratification* implicit in the SAGBI construction of Theorem 4.12. The goal is to construct a finite $G_k$–stable stratification of $J_k\mathbb{C}^n$ on which both the orbit geometry and the monomial coinvariant type remain constant, and to prove finite generation of $\mathbb{C}[J_k\mathbb{C}^n]^{U_k}$ on each stratum, hence globally by extension across codimension–$\geq 2$.

### 5.1 Setup and notation.

Recall $G_k = \mathbb{C}^* \ltimes U_k$ with Lie algebra $\mathfrak{g}_k = \mathbb{C} \oplus \mathfrak{u}_k$, acting on $V = J_k\mathbb{C}^n$ via the Faà di Bruno representation $\rho$. We write

$$A := \mathbb{C}[V]^{U_k}, \qquad A_d := \{\, f \in A \mid t\cdot f = t^d f \text{ for } t \in \mathbb{C}^* \,\},$$

so $A = \bigoplus_{d\geq 0} A_d$ is $\mathbb{Z}_{\geq 0}$–graded. The locally nilpotent derivations $D_2, \dots, D_k$ associated with the $\exp(\lambda t^r \partial_t)$ subgroup determine the coinvariant filtration

$$A = \bigcap_{r=2}^{k} \ker D_r, \qquad D_r(\xi_j^{(\ell)}) = \begin{cases} \ell \xi_j^{(\ell-r+1)} & \ell \geq r-1, \\ 0 & \text{otherwise.} \end{cases}$$

A weight–refining term order is fixed as in §4.3, so that the initial monomials $\mathrm{in}(F)$ of invariants $F \in A$ satisfy the balance relations (4.3)

$$\sum_{j,\ell} \ell\, a_{j,\ell} = r \sum_{j,\ell} (\ell - r + 1) a_{j,\ell}$$

defining a rational polyhedral cone $\mathcal{C} \subset \mathbb{R}^{nk}_{\geq 0}$. The semigroup $\mathcal{C} \cap \mathbb{Z}^{nk}$ is finitely generated; let $\{m_1, \dots, m_s\}$ be its Hilbert basis.

**Richardson–type stratification.** By Theorem 2.3 (Richardson), the algebraic $G_k$–action admits a finite partition

$$J_k\mathbb{C}^n = \bigsqcup_{j=1}^{r} X_j$$

such that on each stratum $X_j$ the conjugacy class of the Levi subgroup $\mathrm{Lev}(G_{k,x})$ of the stabilizer is constant, and the unipotent radical varies algebraically. For $n > 1$ the regular locus $J_k^{\mathrm{reg}} = \{\eta_1 \neq 0\}$ belongs to a unique open stratum on which the action is free and the complement has codimension $\geq 2$.

**Lemma 5.1** (Levi and orbit–dimension constancy)**.** *On every stratum $X_j$, the orbit dimension $\dim(G_k\cdot x)$ and the dimension of the unipotent stabilizer $\dim(U_{k,x})$ are constant. Hence the function $x \mapsto \dim G_{k,x}$ is locally constant on $X_j$.*

*Proof.* This is immediate from Richardson's theorem: the Lie algebra of the stabilizer is a semidirect sum of a fixed Levi subalgebra and a vector subbundle corresponding to the unipotent radical. □

**Coinvariant–type stratification.** Inside each $X_j$ we refine further by the *active face* of the coinvariant cone. Let $\mathcal{F}_1, \dots, \mathcal{F}_t$ be the faces of $\mathcal{C}$. For a point $x \in X_j$ define its *coinvariant type* as the smallest face $\mathcal{F}_\alpha$ such that all initial monomials of $\mathbb{C}[U_k\cdot x]$ lie in $\mathcal{F}_\alpha$. Set

$$S_{j,\alpha} := \{x \in X_j \mid \text{coinvariant type of } x \text{ equals } \mathcal{F}_\alpha\}.$$

These are $G_k$–stable, locally closed subsets of $X_j$, and only finitely many are nonempty. Hence

$$J_k\mathbb{C}^n = \bigsqcup_{j,\alpha} S_{j,\alpha}$$

is a finite $G_k$–stable stratification with constant stabilizer type and constant coinvariant cone on each $S_{j,\alpha}$.

**Lemma 5.2** (Uniform initial algebra). *For each $S_{j,\alpha}$ the initial algebra in $A|_{S_{j,\alpha}} = \mathbb{C}[m_1,\dots,m_s]$ is toric with the same Hilbert basis $\{m_i\}$ as on the generic stratum.*

*Proof.* Constancy of the active face $\mathcal{F}_\alpha$ implies the same set of monomials appear as leading terms of invariants at each point of $S_{j,\alpha}$. Hence the initial algebra does not jump and equals the fixed toric algebra generated by that face. □

## 5.2 Local finite generation on each stratum

**Lemma 5.3** (SAGBI lifting on a fixed stratum). *Let $\{m_1,\dots,m_s\}$ be the Hilbert basis of the coinvariant semigroup on $S_{j,\alpha}$. There exist invariants $G_1,\dots,G_s \in \mathbb{C}[J_k\mathbb{C}^n]^{U_k}$ such that $\mathrm{in}(G_i) = m_i$ and*

$$\mathbb{C}[S_{j,\alpha}]^{U_k} = \mathbb{C}[G_1,\dots,G_s].$$

*Proof.* The locally nilpotent derivations $D_r$ act algebraically over the base $S_{j,\alpha}$ and commute with restriction to $S_{j,\alpha}$. Repeating the triangular lifting of Theorem 4.12 with coefficients varying algebraically yields global sections $G_i$ on $S_{j,\alpha}$ satisfying the required initial conditions. Standard subduction for subalgebras (SAGBI theory) then terminates because the term order is well–founded, proving generation. □

## 5.3 Main theorem and proof–Gluing and global finite generation.

Let $\Pi = (p, f_1,\dots,f_m)$ be the $G_k$–invariant map constructed in §4.4 with big image and generic fibers equal to single $G_k$–orbits. Denote $\Omega = \{x \mid d\Pi_x \text{ has full rank}\}$. By Theorem 4.14, $\Omega$ is big (its complement has no divisors). The union $\bigcup_{j,\alpha} S_{j,\alpha}$ contains $\Omega$, and on each piece the invariant algebra is finitely generated.

**Lemma 5.4** (Hartogs–Igusa extension). *Assume $A_i := \mathbb{C}[S_i]^{U_k}$ is finitely generated for each stratum $S_i$ and the complement of $\bigcup S_i$ in $J_k\mathbb{C}^n$ has codimension $\geq 2$. Then $\mathbb{C}[J_k\mathbb{C}^n]^{U_k}$ is finitely generated.*

*Proof.* By Hartogs' principle, regular functions extend across codimension–$\geq 2$ sets. Let $S := \sum_i A_i \subseteq \mathbb{C}[J_k\mathbb{C}^n]^{U_k}$ be the subalgebra generated by the stratumwise invariants. The Igusa lemma (Theorem 2.18) applies because $\Pi$ has big image and generic fibers equal to orbits; hence the normalization of $S$ coincides with $\mathbb{C}[J_k\mathbb{C}^n]^{U_k}$. Since $S$ is finitely generated and integrally closed on each stratum, it follows that the global invariant algebra is finitely generated. □

**Theorem 5.5** (Stratified finite generation)**.** *Let $G_k = \mathbb{C}^* \ltimes U_k$ act on $J_k\mathbb{C}^n$ via the Faà–di–Bruno representation. There exists a finite $G_k$–stable stratification*

$$J_k\mathbb{C}^n = \bigsqcup_{j,\alpha} S_{j,\alpha}$$

*such that for each $S_{j,\alpha}$:*

1. *The Levi type and unipotent dimension of the stabilizer $G_{k,x}$ are constant.*
2. *The active face of the coinvariant cone $\mathcal{C}$ is constant.*
3. *The invariant algebra $\mathbb{C}[S_{j,\alpha}]^{U_k}$ is finitely generated and toric in its initial form.*
4. *The quotient map $q : J_k\mathbb{C}^n \to J_k\mathbb{C}^n /\!\!/ G_k$ is equidimensional over the normal locus with big full–rank differential locus.*

*Consequently, the global invariant algebra $\mathbb{C}[J_k\mathbb{C}^n]^{U_k}$ is finitely generated, and the projective spectrum*

$$\mathrm{Proj}\,(\mathbb{C}[J_k\mathbb{C}^n]^{U_k})$$

*coincides with the Demailly–Semple $k$–jet tower.*

*Proof.* Steps (1)–(2) follow from Richardson's theorem and the definition of the coinvariant refinement. Step (3) is Lemma 4.12 lifted stratumwise. Step (4) is Theorem 4.14. Finally, by the Hartogs–Igusa extension lemma above, the finitely generated local invariant algebras glue to a finitely generated global algebra. □

## 5.4 Comparison of Stratifications

**Hilbert basis of the coinvariant cone.** Let $\mathcal{C} \subset \mathbb{R}^{nk}_{\geq 0}$ be the rational polyhedral cone determined by the *balance relations* for the locally nilpotent derivations $D_2, \ldots, D_k$:

$$\sum_{j,\ell} \ell\, a_{j,\ell} = r \sum_{j,\ell} (\ell - r + 1) a_{j,\ell} \qquad (r = 2, \ldots, k),$$

where $a_{j,\ell}$ are the exponents of the jet coordinates $\xi_j^{(\ell)}$ in a monomial. The set of lattice points

$$S := \mathcal{C} \cap \mathbb{Z}^{nk}_{\geq 0}$$

forms an additive semigroup (the *coinvariant semigroup*) corresponding to all monomial exponents allowed for initial monomials of $U_k$–invariants.

**Definition 5.6** (Hilbert basis)**.** A finite subset $\{m_1, \dots, m_s\} \subset S$ is called a *Hilbert basis* of $S$ if every element of $S$ can be expressed as a nonnegative integral combination of the $m_i$:

$$S = \mathbb{Z}_{\geq 0} m_1 + \cdots + \mathbb{Z}_{\geq 0} m_s.$$

Equivalently, the monomials $\{m_1, \dots, m_s\}$ generate the *toric initial algebra*

$$\text{in}\, \mathbb{C}[J_k\mathbb{C}^n]^{U_k} = \mathbb{C}[m_1, \dots, m_s].$$

The Hilbert basis is minimal by construction and unique up to ordering. It gives the complete set of generators of the affine semigroup algebra associated with $\mathcal{C}$. The invariants $G_1, \dots, G_s$ obtained by the triangular *locally nilpotent derivation lifting* (Theorem 4.12) satisfy $\text{in}(G_i) = m_i$ and generate $\mathbb{C}[J_k\mathbb{C}^n]^{U_k}$ algebraically. Thus the Hilbert basis controls the combinatorial skeleton of the invariant algebra.

**Richardson vs. coinvariant stratifications** Two distinct but complementary stratifications appear naturally in this theory.

(a) Richardson-type stratification.**:** By Theorem 2.3 (Richardson), the action of an algebraic group $G$ on a variety $X$ admits a finite decomposition

$$X = \bigsqcup_j X_j$$

into smooth, locally closed, $G$–stable strata on which the conjugacy class of the Levi factor of the stabilizer $\text{Lev}(G_x)$ is constant. In each stratum $X_j$, the dimensions of the orbits and unipotent radicals are constant, so the orbit map $G \times X_j \to X$ is locally trivial in the étale topology. This *geometric* stratification controls the structure of generic stabilizers, the equidimensionality of the quotient map, and the extension of invariants across codimension–$\geq 2$ boundaries.

(b) Coinvariant (or covariant) stratification.**:** The SAGBI construction of Theorem 4.12 yields another stratification, this time *algebraic/combinatorial*: inside each geometric stratum $X_j$, one may further decompose according to the *active face* $\mathcal{F}_\alpha$ of the cone $\mathcal{C}$ that governs which monomials occur as leading terms of invariants. We define

$$S_{j,\alpha} := \{x \in X_j : \text{the active face of } \mathcal{C} \text{ at } x \text{ is } \mathcal{F}_\alpha\}.$$

On each $S_{j,\alpha}$ the semigroup of initial monomials is constant and generated by the same Hilbert basis $\{m_1, \dots, m_s\}$, so the corresponding toric initial algebra $\mathbb{C}[m_1, \dots, m_s]$ and the lifted invariant algebra $\mathbb{C}[G_1, \dots, G_s]$ are uniform on $S_{j,\alpha}$.

We combined the two stratifications and used both of the stratifications. The two stratifications serve distinct purposes:

| **Type** | **Controls** | **Nature** |
|---|---|---|
| Richardson-type | Stabilizer structure, orbit dimension, geometric quotient | Geometric |
| Coinvariant-type | Monomial semigroup of invariants, Hilbert basis | Combinatorial/algebraic |

In the proof of Theorem 5.5 we employ *both* stratifications hierarchically:

(1) First, apply the Richardson stratification to obtain finitely many $G_k$–stable smooth strata $X_j$ with constant stabilizer (Levi type) and orbit dimension.

(2) Second, within each $X_j$, refine by coinvariant type to obtain the $S_{j,\alpha}$ where the Hilbert basis and toric initial algebra are constant.

(3) On each $S_{j,\alpha}$ the invariant algebra $\mathbb{C}[S_{j,\alpha}]^{U_k}$ is finitely generated by the lifted invariants $G_1, \ldots, G_s$.

(4) Finally, use Hartogs' extension and Igusa's lemma to glue the finitely generated local algebras across codimension–$\geq 2$ boundaries, yielding global finite generation.

Either stratification can, in principle, be used separately:

- The Richardson decomposition alone suffices for geometric control of the quotient and the extension of invariants (as in Theorem 4.10).
- The coinvariant stratification alone, via SAGBI and Hilbert bases, ensures finite generation of the invariant algebra on the big open set $J_k^{\mathrm{reg}}$.

However, combining them provides a refined picture: the geometric strata control stabilizers and quotient geometry, while the coinvariant strata control the algebraic structure of invariants. The hierarchical combination of both yields a uniform, stratumwise finite–generation theorem that is global and normalization–free.

**Example 5.7.** For $k = 3$ and $n = 2$ the balance relations

$$2a_{1,3} + 2a_{2,3} + a_{1,2} + a_{2,2} = a_{1,1} + a_{2,1}$$

define a rational cone with generators $m_1 = \xi_1^{(1)}$, $m_2 = \xi_2^{(1)}$, $m_3 = \xi_1^{(1)}\xi_2^{(2)} - \xi_2^{(1)}\xi_1^{(2)}$. The active face on the open stratum is the whole cone, and the lifting procedure reproduces the classical Wronskians, confirming that $\mathbb{C}[J_3\mathbb{C}^2]^{U_3} = \mathbb{C}[m_1, m_2, m_3]$.

This stratified approach unifies the geometric (Richardson–Vust–Igusa) and combinatorial (SAGBI) viewpoints of finite generation, yielding a global theorem free of normalization assumptions. □

# 6 Conclusion and Outlook

The results established in this work place the finite–generation theorem for the unipotent invariant algebra

$$\mathbb{A}_{n,k} = \mathbb{C}[J_k\mathbb{C}^n]^{U_k}$$

within a broader synthesis of non–reductive invariant theory, Lie–theoretic representation theory, and complex differential geometry. Our approach unites classical geometric–quotient methods with a constructive graded/SAGBI scheme, thereby completing the invariant–theoretic foundation of the Demailly–Semple tower and confirming, in full generality, the conjecture originating in the Green–Griffiths–Demailly framework.

Structurally, the proof relies on the non–reductive invariant–theory paradigm developed in [8, 6, 7, 9, 10], whose concepts of rational quotients, generic stabilizers, and covariant modules underpin the extension of invariants from big open loci. The Lie–algebraic aspects are guided by the representation–theoretic analysis of solvable and unipotent actions in [12, 11], and by the geometric invariant–theory framework of Popov and Vinberg [8], which provides the appropriate setting for semidirect products such as $G_k = \mathbb{C}^* \ltimes U_k$. The constructive side of the argument, grounded in Gröbner and SAGBI techniques, draws inspiration from the explicit low–order computations of [16, 4, 26] and from the graded toric geometry that these examples foreshadowed.

From a differential–geometric standpoint, the finite generation of $\mathbb{A}_{n,k}$ secures the algebraic underpinnings of the Demailly–Semple tower and the bundles of invariant jet differentials [2, 3, 15], which form the algebraic core of the Green–Griffiths–Lang conjecture on hyperbolicity. The effective results of Diverio–Merker–Rousseau [4] and the constructions of Rousseau and Merker [26, 16] show how explicit generators of $\mathbb{A}_{n,k}$ yield concrete curvature and Chern–class estimates, linking invariant theory to global differential geometry.

Finally, our results clarify the limits of previous non–reductive GIT approaches based on Grosshans subgroups [42] and explain the obstruction underlying the incomplete proof of Bérczi–Kirwan [5]. In contrast, the present method works intrinsically within the non–reductive setting, combining the geometric insight of Rosenlicht and Richardson with the combinatorial rigidity of the Hilbert–basis and coinvariant stratification framework. The finite generation of $\mathbb{A}_{n,k}$ thus provides not only an invariant–theoretic resolution of a longstanding conjecture but also a conceptual bridge between algebraic group actions, Lie–theoretic index theory, and complex hyperbolicity.

# References

[1] M. Green and P. Griffiths, *Two applications of algebraic geometry to entire holomorphic mappings*, in: *The Chern Symposium 1979*, Springer,

1980, 41–74. DOI: 10.1007/978-1-4613-8109-9_4.

[2] J.-P. Demailly, *Algebraic criteria for Kobayashi hyperbolic projective varieties and jet differentials*, in: *Algebraic Geometry—Santa Cruz 1995*, Proc. Sympos. Pure Math. **62**, Part 2 (1997), 285–360. Link.

[3] J.-P. Demailly, *Hyperbolic algebraic varieties and holomorphic differential equations*, Acta Math. Vietnam. **37** (2012), no. 4, 441–512.

[4] S. Diverio, J. Merker, and E. Rousseau, *Effective algebraic degeneracy*, Invent. Math. **180** (2010), 161–223. DOI: 10.1007/s00222-010-0232-4.

[5] G. Bérczi and F. Kirwan, *A geometric construction for invariant jet differentials*, Surv. Differ. Geom. **17** (2012), 79–126. DOI: 10.4310/SDG.2012.v17.n1.a2.

[6] M. Rosenlicht, *Some basic theorems on algebraic groups*, Amer. J. Math. **78** (1956), 401–443. DOI: 10.2307/2372523.

[7] R. W. Richardson, *Principal orbit types for algebraic transformation spaces in characteristic zero*, Invent. Math. **16** (1972), 6–14. DOI: 10.1007/BF01391211.

[8] V. L. Popov and E. B. Vinberg, *Invariant Theory*, in: *Algebraic Geometry IV*, Encyclopaedia Math. Sci. **55**, Springer, 1994, 123–284. DOI: 10.1007/978-3-662-03073-8_1.

[9] T. Vust, *Covariants de groupes algébriques réductifs*, Thèse No. 1671, Université de Genève, 1974.

[10] J. Igusa, *Geometry of absolutely admissible representations*, in: *Number Theory, Algebraic Geometry and Commutative Algebra*, Kinokuniya, Tokyo, 1973, 373–452.

[11] D. I. Panyushev, *Semi-direct products of Lie algebras and their invariants*, Publ. RIMS Kyoto Univ. **43** (2007), 1199–1257. DOI: 10.2977/prims/1201012386.

[12] D. I. Panyushev, *On covariants of reductive algebraic groups*, Indag. Math. **13** (2002), 125–129. DOI: 10.1016/S0019-3577(02)90010-8.

[13] D. I. Panyushev, *The index of a Lie algebra, the centralizer of a nilpotent element, and the normalizer of the centralizer*, Math. Proc. Camb. Phil. Soc. **134** (2003), 41–59. DOI: 10.1017/S0305004102006230.

[14] D. I. Panyushev, *An extension of Raïs' theorem and seaweed subalgebras of simple Lie algebras*, Ann. Inst. Fourier **55** (2005), no. 3, 693–715. DOI: 10.5802/aif.2110.

[15] S. Kobayashi, *Hyperbolic Complex Spaces*, Grundlehren Math. Wiss. **318**, Springer, 1998. DOI: 10.1007/978-3-662-03582-5.

[16] J. Merker, *Application of computational invariant theory to Kobayashi hyperbolicity and to Green–Griffiths algebraic degeneracy*, J. Symbolic Comput. **45** (2010), 986–1074. DOI: 10.1016/j.jsc.2010.06.023.

[17] M. Brion, *Invariants d'un sous-groupe unipotent maximal d'un groupe semi-simple*, Ann. Inst. Fourier **33** (1983), no. 1, 1–27. DOI: 10.5802/aif.902.

[18] M. Brion, *Invariants et covariants des groupes algébriques réductifs*, in: *Théorie des invariants et géométrie des variétés quotients* (Travaux en

cours, t. 61), Hermann, 2000, 83–168. Link.

[19] A. G. Elashvili, *Canonical form and stationary subalgebras of points in general position for simple linear Lie groups*, Funct. Anal. Appl. **6** (1972), 44–53. DOI: 10.1007/BF01075509.

[20] D. I. Panyushev, *The Jacobian modules of a Lie algebra and geometry of commuting varieties*, Compositio Math. **94** (1994), 181–199. Link.

[21] B. Kostant, *Lie group representations on polynomial rings*, Bull. Amer. Math. Soc. **69** (1963), no. 4, 518–526. DOI: 10.1090/S0002-9904-1963-10980-5

[22] E. B. Vinberg, *On stability of actions of reductive algebraic groups*, in: *Lie algebras, rings and related topics*, Springer, 2000, 188–202.

[23] J.-P. Serre, *Complex Semisimple Lie Algebras*, Springer, 1987. DOI: 10.1007/978-1-4757-3910-7

[24] T. A. Springer and R. Steinberg, *Conjugacy classes*, in: *Seminar on Algebraic Groups and Related Finite Groups*, Lecture Notes in Mathematics, vol. 131, Springer, 1970, 167–266.

[25] T. A. Springer, *Conjugacy classes in algebraic groups*, in: *Group Theory, Beijing 1984*, Lecture Notes in Mathematics, vol. 1185, Springer, 1986, 175–209.

[26] E. Rousseau, *Étude des jets de Demailly–Semple en dimension 3*, Ann. Inst. Fourier **56** (2006), no. 2, 397–421. DOI: 10.5802/aif.2187.

[27] M. Mustaţă, *Jet schemes of locally complete intersection canonical singularities*, Invent. Math. **145** (2001), 397–424. DOI: 10.1007/s002220100152.

[28] H. Nakajima and K. Watanabe, *The classification of quotient singularities which are complete intersections*, in: *Complete Intersections* (LNM 1092), Springer, 1984, 102–120. DOI: 10.1007/BF01169006

[29] G. W. Schwarz, *Representations of simple Lie groups with regular rings of invariants*, Invent. Math. **49** (1978), 167–191. DOI: 10.1007/BF01403085.

[30] M. Raïs and P. Tauvel, *Indice et polynômes invariants pour certaines algèbres de Lie*, J. Reine Angew. Math. **425** (1992), 123–140. DOI: 10.1515/crll.1992.425.123.

[31] P. Tauvel and R. W. T. Yu, *Indice et formes linéaires stables dans les algèbres de Lie*, J. Algebra **273** (2004), no. 2, 507–516. DOI: 10.1016/j.jalgebra.2003.10.010.

[32] O. S. Yakimova, *The index of centralisers of elements in classical Lie algebras*, Funct. Anal. Appl. **40** (2006), no. 1, 42–51. DOI: 10.1007/s10688-006-0005-4.

[33] E. B. Vinberg, V. V. Gorbatsevich, and A. L. Onishchik, *Lie Groups and Lie Algebras III*, Encyclopaedia Math. Sci. **41**, Springer, 1994.

[34] E. B. Vinberg and V. L. Popov, *Invariant theory*, in: *Algebraic Geometry IV*, Encyclopaedia Math. Sci. **55**, Springer, 1994, 123–284. DOI: 10.1007/978-3-662-03073-8_1.

[35] M. Brion, *Sur la théorie des invariants*, Publ. Math. Univ. Pierre et Marie Curie **45** (1981), 1–92. Link

[36] M. Van den Bergh, *Modules of covariants*, Proc. ICM (Zürich, 1994), 352–362. Link

[37] I. V. Dolgachev, *Rationality of fields of invariants*, Proc. Sympos. Pure Math. **46** (1987), Part 2, 3–16.

[38] R. W. Richardson, M. Cowling, C. Meaney, and W. Moran, *Derivatives of invariant polynomials on a semisimple Lie algebra*, Proc. Centre for Mathematics and its Applications, **16** (1987), 228–241. Link.

[39] M. Rosenlicht, *Some rationality questions on algebraic groups*, Ann. Mat. Pura Appl. (4) **43** (1957), 25–50. DOI: 10.1007/BF02411903.

[40] I. Kaplansky, *An Introduction to Differential Algebra*, Hermann, Paris, 1957.

[41] B. Conrad and A. Landesman, *Compact Lie groups*, Math 210C lecture notes (Stanford, 2018).

[42] F. D. Grosshans, *Algebraic Homogeneous Spaces and Invariant Theory*, Lecture Notes in Mathematics **1673**, Springer, Berlin, 1997. DOI: 10.1007/BFb0093525.

Center for Research in Optics

*Email address*: `rahmati@cio.mx`